\documentclass[preprint,12pt]{elsarticle}

\usepackage{amssymb}
\usepackage{amsmath}

\usepackage{comment}
\usepackage{caption}
\usepackage{url}

\DeclareMathOperator{\trace}{trace}

\usepackage{booktabs} 
\usepackage{siunitx}  

\newcommand{\bfA}[0]{\boldsymbol{A}}

\newcommand{\bfe}[0]{\boldsymbol{e}}

\newcommand{\bfg}[0]{\boldsymbol{g}}

\newcommand{\bfI}[0]{\boldsymbol{I}}

\newcommand{\bfJ}[0]{\boldsymbol{J}}

\newcommand{\bfn}[0]{\boldsymbol{n}}

\newcommand{\bfp}[0]{\boldsymbol{p}}

\newcommand{\bfu}[0]{{\boldsymbol{u}}}
\newcommand{\bfU}[0]{\boldsymbol{U}}
\newcommand{\bfv}[0]{\boldsymbol{v}}

\newcommand{\bfx}[0]{\boldsymbol{x}}

\newcommand{\sigt}[0]{\tilde{\sigma}}

\newcommand{\hatn}[0]{\hat{\bfn}}

\newcommand{\ccal}[0]{\mathcal{C}}
\newcommand{\ecal}[0]{\mathcal{E}}
\newcommand{\jcal}[0]{\mathcal{J}}
\newcommand{\kcal}[0]{\mathcal{K}}
\newcommand{\mcal}[0]{\mathcal{M}}
\newcommand{\ncal}[0]{\mathcal{N}}
\newcommand{\ocal}[0]{\mathcal{O}}
\newcommand{\pcal}[0]{\mathcal{P}}
\newcommand{\qcal}[0]{\mathcal{Q}}
\newcommand{\scal}[0]{\mathcal{S}}
\newcommand{\tcal}[0]{\mathcal{T}}

\newcommand{\vcal}[0]{\mathcal{V}}

\newcommand{\ddt}[1]{\frac{\partial #1}{\partial t}}

\newcommand{\visc}[0]{{\eta}}
\newcommand{\phitest}[0]{\varphi_\phi}
\newcommand{\mutest}[0]{\varphi_\mu}

\newcommand{\shape}[0]{\varphi}
\newcommand{\vshape}[0]{\boldsymbol{\varphi}}

\newcommand{\intom}[1]{\int_\Omega{#1}\,d\bfx}

\newcommand{\NL}[0]{\text{NL}}

\newcommand{\nst}[0]{N_{st}}
\newcommand{\ncalst}[0]{\ncal_{st}}
\newcommand{\navg}[0]{N_{avg}}
\newcommand{\ncalavg}[0]{\ncal_{avg}}
\newcommand{\nI}[0]{{n_I}}
\newcommand{\nt}[0]{{n_T}}

\usepackage{color}
\usepackage[dvipsnames,svgnames]{xcolor}
\newcommand{\inred}[1]{{\color{red} #1}}
\newcommand{\ingreen}[1]{{\color{DarkGreen} #1}}
\newcommand{\inblue}[1]{{\color{blue} #1}}

\usepackage{cleveref}
\crefformat{equation}{(#2#1#3)}
\crefrangeformat{equation}{(#3#1#4) to~(#5#2#6)}
\crefmultiformat{equation}{(#2#1#3)}%
{ and~(#2#1#3)}{, (#2#1#3)}{ and~(#2#1#3)}

\journal{Journal of Computational Physics}

\begin{document}

\begin{frontmatter}


\title{Anisotropic mesh adaptation for unsteady two-phase flow simulation
with the Cahn-Hilliard Navier-Stokes model\tnoteref{label1}}
\author[1]{A. Bawin\corref{cor1}}
\ead{arthur.bawin@polymtl.ca}
\author[1]{S. Étienne}
\ead{stephane.etienne@polymtl.ca}
\author[1]{C. Béguin}
\ead{c.beguin@polymtl.ca}
\cortext[cor1]{Corresponding author}

\affiliation[1]{organization={Département de génie mécanique,
Polytechnique Montréal},
addressline={2500 Chem. de Polytechnique},
city={Montreal},
citysep={}, 
%
postcode={H3T 0A3},
state={QC},
country={Canada}}

\begin{abstract}
We present an anisotropic mesh adaptation procedure based on Riemannian metrics for the simulation of two-phase incompressible flows with non-matching densities.
The system dynamics are governed by the Cahn-Hilliard Navier-Stokes (CHNS) equations,
discretized with mixed finite elements and implicit time-stepping.
Spatial accuracy is controlled throughout the simulation by the \emph{global transient fixed-point method} from Alauzet \emph{et al.},
in which the simulation time is divided into sub-intervals, each associated with an adapted anisotropic mesh.
The simulation is run in a fixed-point loop until convergence of each mesh--solution pair.
Each iteration takes advantage of the previously computed solution and accurately predicts the flow variations.
This ensures that the mesh always captures the fluid-fluid interface,
and allows for a dynamic control of the interface thickness at a fraction of the computational cost compared to uniform or isotropic grids.
Moreover, using a modest number of time sub-intervals reduces the transfer error from one mesh to another,
which would otherwise eventually spoil the numerical solution.
The overall adaptive procedure is verified with manufactured solutions and the well-known rising bubble benchmark.
\end{abstract}




\end{frontmatter}


\section{Introduction}
\label{sec:intro}
Two-phase flows are involved in a variety of natural phenomena, such as
ocean-atmosphere interactions or volcanic eruptions, as well as industrial applications,
such as microfluidics, cooling systems or liquefied gas transport.
In recent years, phase-field methods, also called diffuse interface methods,
have proven flexible and reliable for the numerical simulation of two- and multiphase flows.
These Eulerian methods model the fluid-fluid interface by a smooth transition region of
finite thickness \citep{garcia2025numerical}.
An accurate description of the interfacial region is
essential to capture phase interactions and provide reliable predictions.
Representing the interface remains challenging, however:
using uniform grids for unsteady simulations helps handle the interface
motion without requiring prior knowledge of its evolution, but resolving the
interface thickness imposes a stringent constraint on the mesh size throughout the domain.
This yields costly computations,
with unnecessary degrees of freedom in single-phase regions.
To mitigate this, adaptive mesh refinement (AMR) techniques
were proposed in, e.g., \citep{hintermuller2013adaptive, chen2016efficient, khanwale2020simulating, alphonius2025lethe}.
In these works, spatial discretization is handled with simplicial or tree-based meshes,
and AMR is performed with a mark-and-refine approach based on various a posteriori error estimators.
These methods yield adapted isotropic meshes, in which the size of the elements
vary, but where each element exhibits essentially the same size in all directions.

Further reduction of degrees of freedom can be achieved by
considering anisotropic meshes, characterized by stretched elements along
physics-induced directions of anisotropy.
Anisotropic mesh adaptation has proven to be computationally efficient to capture
flow fields with steep gradients (e.g., incompressible flows \citep{coupez2013solution, bawin2024metric} and multiphase flows \citep{alauzet2018time})
or discontinuities (e.g., inviscid Euler equations \citep{loseille2008adaptation, alauzet2016decade}).
Anisotropic mesh adaptation is particularly appealing for multiphase flows,
as the interface is an obvious target for refinement.
To our best knowledge, however, the only examples of anisotropic mesh adaptation for
unsteady simulation with phase-field models,
and particularly with a Cahn-Hilliard Navier-Stokes (CHNS) model,
are presented in \citep{pigeonneau2019discontinuous, bonart2021simulation}.
A shortcoming of \citep{pigeonneau2019discontinuous}
is that 
the mesh is refined at each time step.
While this ensures that the mesh does not lag behind the solution,
remeshing every single step prevents reusing the solver structures,
which can be costly to compute (e.g.,
preconditioner or linearized Jacobian matrix for the nonlinear solver).
Moreover, the solution must be interpolated at each time step,
leading to the accumulation of transfer error.

To mitigate these issues, we consider the \emph{global transient fixed-point}
algorithm presented in \citep{alauzet2018time},
which relies on converging the mesh--solution pair in a fixed-point method.
A transient solution is computed on an initial (e.g., uniform) grid,
then the mesh is adapted at equally spaced time intervals (typically 10 to 20),
knowing the behavior of the flow.
The number of remeshings is limited to the number of time intervals,
which is typically small compared to the total number of time steps,
and an arbitrary number of time steps can be performed on each adapted mesh.
This approach relies on Riemannian metrics to represent an anisotropic mesh,
and follows the continuous mesh framework introduced in \citep{loseille2011continuous}.
In particular, the adapted mesh on each time interval is \emph{optimal},
in the sense that it minimizes the interpolation error over a chosen sensor
in $L^p(\Omega)$ norm, guaranteeing a global error control.
This method was applied to compressible flows with shocks,
as well as inviscid multiphase (free surface) simulations modeled by the Euler equations \citep{alauzet2018time}.

The goal of this paper is to present a space adaptive procedure for the simulation of
unsteady multiphase flows with a CHNS model written in a volume-averaged formulation.
More precisely,
we present a fully-coupled continuous Galerkin finite element scheme for the CHNS model,
with implicit time stepping (BDF2).
Taylor-Hood and linear elements are used to discretize the velocity-pressure pair,
and the phase marker-chemical potential pair respectively.
Then, we apply the global fixed-point algorithm to obtain
a Hessian-based time-dependent Riemannian metric and perform anisotropic mesh adaptation
within a fixed-point loop.
An emphasis is put on code verification,
through the use of manufactured solutions for both the CHNS flow solver and the mesh adaptive procedure.
We then apply this procedure to the rising bubble benchmark from Hysing et al.
\citep{hysing2009quantitative},
and show that thin and well-resolved interfaces can be captured for a reduced
computational cost.

This paper is organized as follows.
In \Cref{sec:equations}, we describe the Cahn-Hilliard Navier-Stokes system used
in this work for incompressible flows with non-matching densities.
In \Cref{sec:numericalScheme}, the weak formulation of the CHNS system is derived,
followed by the monolithic finite element formulation.
Fundamentals of anisotropic mesh adaptation are recalled in \Cref{sec:anisotropic_adaptation},
while the global transient-fixed point algorithm from \citep{alauzet2018time},
and in particular the time-varying metric used in the applications, is
presented in \Cref{sec:global_fp}.
Verification of both the CHNS solver and the mesh adaptive procedure
are presented in \Cref{sec:numericalApplications} and are followed by numerical
tests of the rising bubble benchmark.

\section{Governing equations}
\label{sec:equations}
We consider a bounded region $\Omega \subset\mathbb{R}^d, d \leq 3,$ filled with two immiscible
and incompressible fluids of densities $\rho_1, \rho_2$ and dynamic viscosities $\visc_1, \visc_2$.
The Cahn-Hilliard Navier-Stokes system
models the dynamic of the mixture through a \emph{phase marker} or \emph{phase field},
a smooth variable $\phi(\bfx,t)$ characterized by
\begin{equation}
  \phi(\bfx,t) =
  \left\{
  \begin{aligned}
      1 ~~~&\text{in fluid 1},\\
      -1 ~~~&\text{in fluid 2}.
  \end{aligned}
\right.
\end{equation}
The transition from fluid 1 to fluid 2 operates in a region of thickness $\ocal(\epsilon)$
with $\epsilon \ll 1$, and the interface can be described as the level set $\phi = 0$.
In the mixture, the density and viscosity are defined as
\begin{equation}
    \rho(\phi) = \frac{\rho_1 - \rho_2}{2} \phi + \frac{\rho_1+\rho_2}{2}, ~~~~~\visc(\phi) = \frac{\visc_1 - \visc_2}{2} \phi + \frac{\visc_1+\visc_2}{2},\\
\end{equation}
i.e., the phase marker can be seen as a volume fraction.

The dynamics of the mixture can be modeled by a combination of the Cahn-Hilliard equation,
describing the dynamics between the two phases,
and the Navier-Stokes equations,
describing the dynamics of the mixture, viewed as a continuum with macroscopic velocity $\bfu$ and pressure $p$.
The Cahn-Hilliard equation
is derived by introducing the free energy of the mixture:
\begin{equation}
  E_\text{mix} = \int_{\Omega} \Psi(\phi,\nabla \phi) \,d\bfx, ~\text{ with }~ \Psi(\phi,\nabla \phi) =  \sigt\epsilon^{-1} F(\phi) + \frac{\sigt\epsilon}{2}\vert \nabla \phi \vert^2,
\end{equation}
where $\sigt = 3\sigma/(2\sqrt{2})$ and $\sigma$ is the fluid-fluid surface tension parameter.
The total energy $E_\text{mix}$ consists of an interfacial energy term proportional to $\vert \nabla \phi \vert^2$,
and an Helmholtz free energy term $F(\phi)$, often approximated by the \emph{double-well} function $F(\phi) = (1 - \phi^2)^2/4$.
We also introduce the chemical potential $\mu$ as the functional derivative:
\begin{equation}
  \begin{aligned}
  \mu = \frac{\delta E_\text{mix}}{\delta\phi} = \frac{\partial \Psi}{\partial \phi} - \nabla \cdot \left( \frac{\partial \Psi}{\partial \nabla \phi} \right)
  &=  \sigt\epsilon^{-1} (\phi^3-\phi) - \sigt\epsilon \Delta \phi.
  \end{aligned}
\end{equation}
%

A variety of formulations are available in the literature for the complete
Cahn-Hilliard Navier-Stokes (CHNS) system, see for instance \citep{ten2023unified}
for a review and unification of some of these models for (quasi-)incompressible flows.
At this point, we should emphasize that the goal of this paper is to present an adaptive
method for phase-field methods
that is solver-agnostic, and can thus be readily applied to all CHNS models.
Consequently, the focus is not on selecting or comparing specific CHNS formulations.
To demonstrate the approach,
we considered the model from Abels et al. \citep{abels2012thermodynamically}:
\begin{align}
  \rho(\phi)\left(\partial_t\bfu + (\bfu \cdot \nabla)\bfu - \bfg\right) + (\bfJ(\mu) \cdot \nabla) \bfu ~~~~~~~~~~~~~~~~~~~~~\,&\nonumber\\[0.2em]
       + \nabla p - \nabla \cdot(2\visc(\phi)d(\bfu)) + \phi\nabla\mu &= 0,
     \label{eq:CHNS_u}\\[0.2em]
     \nabla \cdot \bfu &= 0,\label{eq:CHNS_div}\\[0.2em]
  \partial_t\phi + \bfu \cdot \nabla\phi - \nabla \cdot (M(\phi) \nabla \mu)&= 0,
  \label{eq:CHNS_phi}\\[0.2em]
  \mu -  \sigt\epsilon^{-1} (\phi^3 - \phi) +  \sigt\epsilon \Delta \phi &= 0,
  \label{eq:CHNS_mu}
\end{align}
where $\bfJ(\mu) = M(\phi)(\rho_2-\rho_1)\nabla \mu/2$ is a diffusive flux,
$M(\phi)$ is the isotropic mobility parameter, which is here assumed constant and given by $M = 0.1 \epsilon^2$, similarly to \citep{brunk2025simple}, and with the rate of strain tensor
$d(\bfu) = (\nabla \bfu + (\nabla \bfu)^T)/2$ and gravity acceleration $\bfg = -g\bfe_y$.
Equations \cref{eq:CHNS_phi,eq:CHNS_mu} are the coupled Cahn-Hilliard
gradient flow equation,
whereas \cref{eq:CHNS_u,eq:CHNS_div} are the incompressible Navier-Stokes equations,
with additional terms to account for mass diffusion at the interface
and capillary forces.
This system
is a slightly modified version of the model from \citep{abels2012thermodynamically},
with the Korteweg force $\sigt\epsilon \, \nabla \cdot (\nabla \phi \otimes \nabla \phi)$ replaced by the term
proposed in \citep{ten2023unified}:
\begin{equation}
  \begin{aligned}
  \nabla \cdot \left(\nabla \phi \otimes \frac{\partial\Psi}{\partial\nabla\phi} + (\mu\phi - \Psi)\bfI\right)
  &= \sigt\epsilon \, \nabla \cdot (\nabla \phi \otimes \nabla \phi) + \nabla (\mu\phi - \Psi)\\
  &= \phi \nabla \mu,
\end{aligned}
\end{equation}
which effectively amounts to defining the pressure as $p = p_\text{Abels} - (\mu\phi - \Psi).$

The system \cref{eq:CHNS_u,eq:CHNS_div,eq:CHNS_phi,eq:CHNS_mu} is completed with the following initial and boundary conditions on $\bfu(\bfx,t), \phi(\bfx,t)$ and $\mu(\bfx,t)$:
\begin{equation}
\begin{alignedat}{4}
  \bfu(\bfx, 0)    &= \bfu_0(\bfx) &&\text{ in }\Omega,             ~~~~~~~~~~~~~~~~~~~~~~\phi(\bfx, 0) &&= \phi_0(\bfx)&&\text{ in }\Omega,\\
  \bfu \cdot \hatn &= 0 &&\text{ on }\partial\Omega_\text{slip}, ~~~~~~~~~~~~~~~~\nabla \phi \cdot \hatn &&= 0           &&\text{ on }\partial\Omega,\\
  \bfu             &= 0 &&\text{ on }\partial\Omega_\text{no slip},~~~~~~~~~~~~~   \nabla \mu \cdot \hatn &&= 0           &&\text{ on }\partial\Omega,
  \end{alignedat}
\end{equation}
with $\partial\Omega_\text{slip} \cup \partial\Omega_\text{no slip} = \partial\Omega$ and $\hatn$ the unit outward normal to the boundary $\partial \Omega$.
The zero flux condition on the phase marker $\phi$ amounts to enforcing that
the fluid-fluid interface meets the boundary at an angle of $\pi/2$ everywhere.
Under these boundary conditions, the system preserves the total phase marker
and dissipates the total energy of the system via (e.g., Theorem 2.1 in \citep{chen2021fully}):
\begin{equation}
  \begin{aligned}
  \frac{d}{dt}\int_\Omega{\phi}\,d\bfx &= 0,\\
  \frac{d}{dt}\int_\Omega E_\text{tot}\,d\bfx &= - \int_\Omega \left(M(\phi) \vert \nabla \mu\vert^2 + \visc\, d(\bfu) : \nabla \bfu\right) \,d\bfx \leq 0,
\end{aligned}
\end{equation}
where
\begin{equation}
  E_\text{tot}(\rho,\bfu,\phi) = \kcal + E_\text{mix} = \int_\Omega \left(\frac{\rho}{2}|\bfu|^2 + \sigt\epsilon^{-1} F(\phi) + \frac{\sigt\epsilon}{2} \vert \nabla \phi\vert^2\right) \,d\bfx,
\end{equation}
with $\kcal$ the kinetic energy of the mixture.

\section{Numerical scheme}
\label{sec:numericalScheme}
To solve \cref{eq:CHNS_u,eq:CHNS_div,eq:CHNS_phi,eq:CHNS_mu},
we use the method of lines and adopt a monolithic finite element formulation in space,
using LBB-stable Taylor-Hood mixed elements, and a second-order
implicit backward differentiation (BDF2) discretization in time.
At each time step, the coupled and nonlinear system is solved with a standard Newton-Raphson method.
The discrete numerical scheme is presented in this section.

\subsection{Weak formulation}
To derive the weak formulation of \cref{eq:CHNS_u,eq:CHNS_div,eq:CHNS_phi,eq:CHNS_mu},
the momentum equation \cref{eq:CHNS_u} is multiplied by vector-valued test functions $\bfv$ and integrated over $\Omega$,
while the remaining three scalar equations are multiplied by scalar test functions denoted by $q, \phitest$ and $\mutest$ respectively.
The pressure gradient and all diffusive terms are integrated by parts,
yielding the problem of finding the
quadruplet
$U \triangleq (\bfu, p, \phi, \mu) \in \vcal \times \qcal \times \scal \times \scal$
such that for all
$(\bfv, q, \phitest, \mutest) \in \vcal \times \qcal \times \scal \times \scal$,
we have:
\begin{equation}
\begin{alignedat}{1}
\NL(U)
&\triangleq \left(\rho\partial_t\bfu, \bfv\right)_\Omega
+ \left(\rho(\bfu \cdot \nabla)\bfu, \bfv\right)_\Omega
- (\rho\bfg,\bfv)_\Omega
+ ((\bfJ(\mu) \cdot \nabla) \bfu, \bfv)_\Omega\\[0.4em]
- &(p,\nabla \cdot \bfv)_\Omega + \left< p, \bfv \cdot \hatn\right>_{\partial \Omega}
+ (2\visc d(\bfu), \nabla \bfv) - \left< 2\visc d(\bfu) \cdot \bfn, \bfv \right>_{\partial \Omega}\\[0.4em]
+ &\left(\phi\nabla\mu, \bfv\right)_\Omega
+ (\nabla \cdot \bfu, q)_\Omega\\[0.4em]
+ &(\phi_t, \phitest)_\Omega
+ (\bfu \cdot \nabla \phi, \phitest)_\Omega
+ (M\nabla \mu, \nabla \phitest)_\Omega
- \left< M\nabla\mu \cdot \hatn , \phitest \right>_{\partial \Omega}\\[0.4em]
+ &(\mu, \mutest)_\Omega
- \sigt\epsilon^{-1}(\phi^3 - \phi, \mutest)_\Omega\\[0.4em]
- &(\sigt\epsilon \nabla \phi, \nabla \mutest)_\Omega
+ \sigt\epsilon \left<\nabla\phi \cdot \hatn , \mutest \right>_{\partial \Omega} = 0.
\end{alignedat}
\end{equation}
Here,
$(f,g)_\Omega$ is a generic notation for the inner product $\int_\Omega f \cdot g\,d\bfx$,
where the dot is either a product of scalar fields, a dot product of vector fields or the double contraction of 2-tensor fields,
and
$\left<f,g\right>_{\partial \Omega} \triangleq \int_{\partial \Omega} f \cdot g \,d\bfx$
similarly denotes a boundary integral.
These integrals are meaningful for $\bfu \in (H^1(\Omega))^d$, $p \in L^2(\Omega)$, and $\phi, \mu \in H^1(\Omega)$,
and we define the function spaces:
\begin{equation}
  \begin{aligned}
    \vcal &= \left\{ \bfv \in (H^1(\Omega))^d ~ | ~ \bfv = 0 \text{ on }\Omega_\text{no slip}, \bfv \cdot \hatn = 0 \text{ on }\Omega_\text{slip}\right\},\\
    \qcal &= L^2_0(\Omega) \triangleq \left\{ q \in L^2(\Omega) ~ | ~ \int_\Omega q\,d\bfx = 0 \right\}, ~~~~~\scal = H^1(\Omega).
  \end{aligned}
\end{equation}
The functional $\NL(U)$ is nonlinear in the vector of unknowns $U$,
and $\NL(U) = 0$ is solved with a standard Newton-Raphson method.

\subsection{Finite element and time discretization}
We now consider a conforming simplicial mesh $\tcal$ of $\Omega$, that is,
a conforming triangulation $\tcal = \cup_{i = 1}^{N_e}K_i$ of $\Omega$,
and denote by $\bfu_h, p_h, \phi_h, \mu_h$ the finite element discretizations on $\tcal$ of the unknown fields.
We consider Lagrange finite elements, and denote by $\pcal_k$ the space of continuous and piecewise polynomial functions of degree $k$:
\begin{equation}
  \pcal_k = \left\{ f \in \ccal^0(\Omega)~|~ f|_K \in \mathbb{P}^k(K), ~\forall K \in \tcal\right\},
\end{equation}
where $\mathbb{P}^k(K)$ is the set of polynomials of degree at most $k$ on $K$,
with a similar definition for the vector-valued velocity.
Taylor-Hood $\pcal_2-\pcal_1$ mixed finite elements are used to discretize the velocity-pressure pair,
whereas linear and equal order $\pcal_1-\pcal_1$ elements are used to discretize the $(\phi,\mu)$ pair.
The finite dimensional spaces of interest are thus:
\begin{equation}
  \vcal_h = \pcal_{2} \cap \vcal,~~~
  \qcal_h = \pcal_{1} \cap \qcal,~~~
  \scal_h = \pcal_{1} \cap \scal.
\end{equation}
Introducing global linear and scalar-valued $\shape_i$ and quadratic and vector-valued $\vshape_i$ shape functions,
the discretized velocity and pressure fields write:
\begin{equation}
  \bfu_h(\bfx,t) = \sum_{i = 0}^{N_{\bfu}} \bfU_i(t) \vshape_i(\bfx),
  ~~~~
  p_h(\bfx,t) = \sum_{i = 0}^{N_p} P_i(t) \shape_i(\bfx)
\end{equation}
and similarly for $\phi_h, \mu_h$,
where the sets $\bfU_i$ and $P_i$ are the associated global degrees of freedom.
The time derivatives $\partial_t\bfu$ and $\partial_t\phi$ are discretized with a standard BDF2 method:
\begin{equation}
  \dot{\bfu}_h \triangleq \frac{\partial \bfu_h}{\partial t} \simeq \sum_{j = 0}^2 c_j \bfu_{h}^{n+1-j}, ~~~
  \dot{\phi}_h \triangleq \frac{\partial \phi_h}{\partial t} \simeq \sum_{j = 0}^2 c_j \phi_{h}^{n+1-j},
\end{equation}
where $c_j = [-3/2, 2, -1/2]/\Delta t$.
All quantities at previous times $n$ and $n-1$ are known,
except for the first time step which requires the undefined values $\bfu_h^{-1}$.
To obtain the fields at step $n + 1 = 1$, a single BDF1 step with $\Delta t/10$ is performed,
then a BDF2 step with variable time step and adjusted coefficients is done with the data at times $0, \Delta t/10$ and $\Delta t$.
At each time step, the discrete solution $U_h^{n+1}$ is obtained by solving
the nonlinear problem $\NL(U_h^{n+1}) = 0$ in a monolithic fashion,
which involves solving the following Newton iteration linear system for the increment $\delta U_h^{n+1}$ until convergence:
\begin{equation}
  \jcal(U_{h,m}^{n+1})_{ij} \, (\delta U_{h,m}^{n+1})_j = -\NL_i(U_{h,m}^{n+1}),
\end{equation}
where $m \geq 0$ is the Newton-Raphson iterations counter,
$\jcal(U_{h,m}^{n+1})$ is the linearization of $\NL$ around $U_{h,m}^{n+1}$ (that is,
the Gâteaux derivative of $\NL$ in the direction $\delta U_{h,m}^{n+1}$),
and where the right-hand side reads:
\begin{equation}
\begin{alignedat}{1}
-&\NL_i(U_{h,m}^{n+1})
=
-\intom{\rho^{n+1}_m \left( \dot{\bfu}_m + (\bfu_m^{n+1} \cdot \nabla)\bfu_m^{n+1} - \bfg\right) \cdot \vshape_i}\\[0.4em]
&-\intom{\left((\bfJ(\mu_m^{n+1}) \cdot \nabla)\bfu_m^{n+1} + \phi_m^{n+1} \nabla \mu_m^{n+1}\right)\cdot \vshape_i}\\[0.4em]
&-\intom{\visc_m^{n+1}(\nabla\bfu_m^{n+1} + (\nabla\bfu_m^{n+1})^T) : \nabla \vshape_i - p_m^{n+1}(\nabla \cdot \vshape_i)}\\[0.4em]
&-\intom{(\nabla \cdot \bfu_m^{n+1}) \,\shape_i}\\[0.4em]
&-\intom{\left( \dot{\phi}_m + \bfu_m^{n+1} \cdot \nabla \phi_m^{n+1}\right)\shape_i + M\nabla \mu_m^{n+1} \cdot \nabla \shape_i}\\[0.4em]
&-\intom{\left( \mu_m^{n+1} - \sigt\epsilon^{-1}((\phi_m^{n+1})^3 - \phi_m^{n+1} \right)\shape_i - \sigt\epsilon \nabla \phi_m^{n+1} \cdot \nabla \shape_i},
\end{alignedat}
\end{equation}
with $\rho^{n+1}_m = \rho(\phi_m^{n+1})$ and similarly for the viscosity and where
the subscript $h$ was dropped to avoid cluttering.

\section{Anisotropic mesh adaptation}
\label{sec:anisotropic_adaptation}
We now turn to the mesh adaptive process, once a solution of the CHNS system has been computed.
Anisotropic mesh adaptation is performed using the framework of Riemannian metrics,
where a smoothly varying field of symmetric and positive-definite (SPD) matrices
encodes the target size and orientation of the mesh elements.
In essence, a Riemannian metric modifies the local measures of lengths and areas, and describes,
at the continuous level, triangulations whose elements are elongated along specified directions of anisotropy.
A metric-based anisotropic mesh generator in turn evaluates all lengths and areas with respect to this metric,
and generates quasi-uniform meshes in the metric space, appearing anisotropic in the Euclidean space.
Quasi-uniform, or \emph{quasi-unit}, simplices are characterized by having all
metric-weighted edge lengths within $[1\sqrt{2},\sqrt{2}]$.
Metric-based adaptation is an intrinsically iterative process,
during which the mesh-solution pair is converged
by running the simulation multiple times until little variation in either is achieved.

We consider in particular the \emph{continuous mesh framework} \citep{loseille2011continuous},
which establishes a duality between a discrete mesh and a continuous Riemannian metric.
Continuous operations, such as differentiation and thus optimization,
are carried out directly on the metric, allowing to derive optimal metric fields
for a given target error functional.
Global or local error functionals can be considered to drive the adaptation process \citep{loseille2008adaptation}:
a typical global indicator is the $L^p$-norm of the interpolation error of a given variable.
Adaptation driven by such indicator is also referred to as \emph{feature-based},
and ensures an overall control of the chosen field over the computational domain.
Local indicators include \emph{goal-oriented} functionals,
e.g., mesh adaptation to control the lift or drag coefficient, and typically require solving an adjoint problem.
In this work, we consider feature-based adaptation, and aim at controlling the
global interpolation error over the phase marker $\phi$.
As a result, the adapted meshes accurately captures the interface,
but the adaptation is mostly blind to the
the flow features in the single-component regions where $\phi = \pm 1$.
Another interesting  choice would be to adapt with respect to the velocity (e.g., $\Vert \bfu \Vert$)
and use the intersected metric $\mcal = \mcal_\phi \cap \mcal_\bfu$ to drive the adaptation.\mbox{}\\

In the rest of this section, we briefly recall the basics of metric-based mesh adaptation,
then present the metric for unsteady simulations used in this paper.
Practical aspects, such as Hessian recovery for $\phi$, size gradation and solution transfer, are also discussed.

\subsection{Basics of metric-based mesh generation}
A Riemannian metric $\mcal$ on a domain $\Omega$ is an application assigning to each point $\bfx$ a symmetric and positive-definite metric tensor, represented by an SPD matrix $\mcal(\bfx) \in \mathbb{R}^{d \times d}$.
At each $\bfx$, the matrix $\mcal(\bfx)$ encodes the local target element sizes $h_i > 0$ and orientations $\bfv_i$ in its eigendecomposition.
In two dimensions, for instance, the matrix representation of a metric writes:
\begin{equation}
  \mcal(\bfx) =
  P\Lambda P^T
  =
  \begin{pmatrix}
    \bfv_1 & \bfv_2
  \end{pmatrix}
  \begin{pmatrix}
    h_1^{-2} & \\
    & h_2^{-2}
  \end{pmatrix}
  \begin{pmatrix}
    \bfv_1 \\ \bfv_2
  \end{pmatrix},
\end{equation}
with $P$ the orthogonal matrix of eigenvectors.
A metric tensor takes two vectors $\bfu,\bfv$ and returns their inner products:
\begin{equation}
  \left<\bfu,\bfv\right>_\mcal = \left<\bfu,\mcal\bfv\right> = \bfu^T\mcal\bfv = \sum_{i,j} u_i\mcal_{ij}v_j,
\end{equation}
where $\left<.,.\right>_\mcal$ and $\left<.,.\right>$ denote the inner product with respect to $\mcal$ and the standard Euclidean metric respectively.
Through the inner product, a metric modifies the computation of lengths as follows.
If $\gamma(t)$ is a curve parameterized by $t \in [0,1]$, then its length with respect to $\mcal$ is the integral:
\begin{equation}
  \ell_\mcal(\gamma) = \int_0^1 \Vert \gamma'(t) \Vert_\mcal\,dt = \int_0^1 \sqrt{\left<\gamma'(t), \gamma'(t)\right>_\mcal}\,dt.
\end{equation}
In particular, a straight mesh edge $\bfe = \bfx_0\bfx_1$ can be parameterized by $\gamma(t) = \bfx_0 + t(\bfx_1-\bfx_0)$,
whose length with respect to $\mcal$ writes:
\begin{equation}
  \ell_\mcal(\bfe) = \int_0^1 \sqrt{(\bfx_1-\bfx_0)^T\mcal(\bfx_0 + t(\bfx_1-\bfx_0))(\bfx_1-\bfx_0)}\,dt.
\end{equation}
This integral can be either evaluated by quadrature,
or simplified by assuming the variation of the metric along the edge \citep{loseille2008adaptation, frey2007mesh}.
The volume of a mesh element $K$ with respect to $\mcal$ is given by 
\begin{equation}
|K|_\mcal = \int_{K} \sqrt{\det \mcal(\bfx)}\,d\bfx. 
\end{equation}
It is closely related to the \emph{mesh complexity} $\ccal$, defined by
\begin{equation}
  \ccal(\mcal) \triangleq \intom{\sqrt{\det \mcal(\bfx)}},
  \label{eq:complexity}
\end{equation}
that is, the total volume of the domain $\Omega$ with respect to $\mcal$.
If $\tcal$ is a triangulation of $\Omega$ that is quasi-unit for $\mcal$ (that is, whose size and orientation are prescribed by $\mcal$),
then the mesh complexity can be viewed as the continuous counterpart to the number of vertices in $\tcal$:
indeed, a large metric determinant $\det \mcal = \lambda_1\lambda_2 = (h_1h_2)^{-2}$ yields small mesh elements in the Euclidean space,
thus a higher count of mesh vertices.

\subsection{Space-time mesh and time-varying metric}
\begin{figure}
  \includegraphics[width=\linewidth]{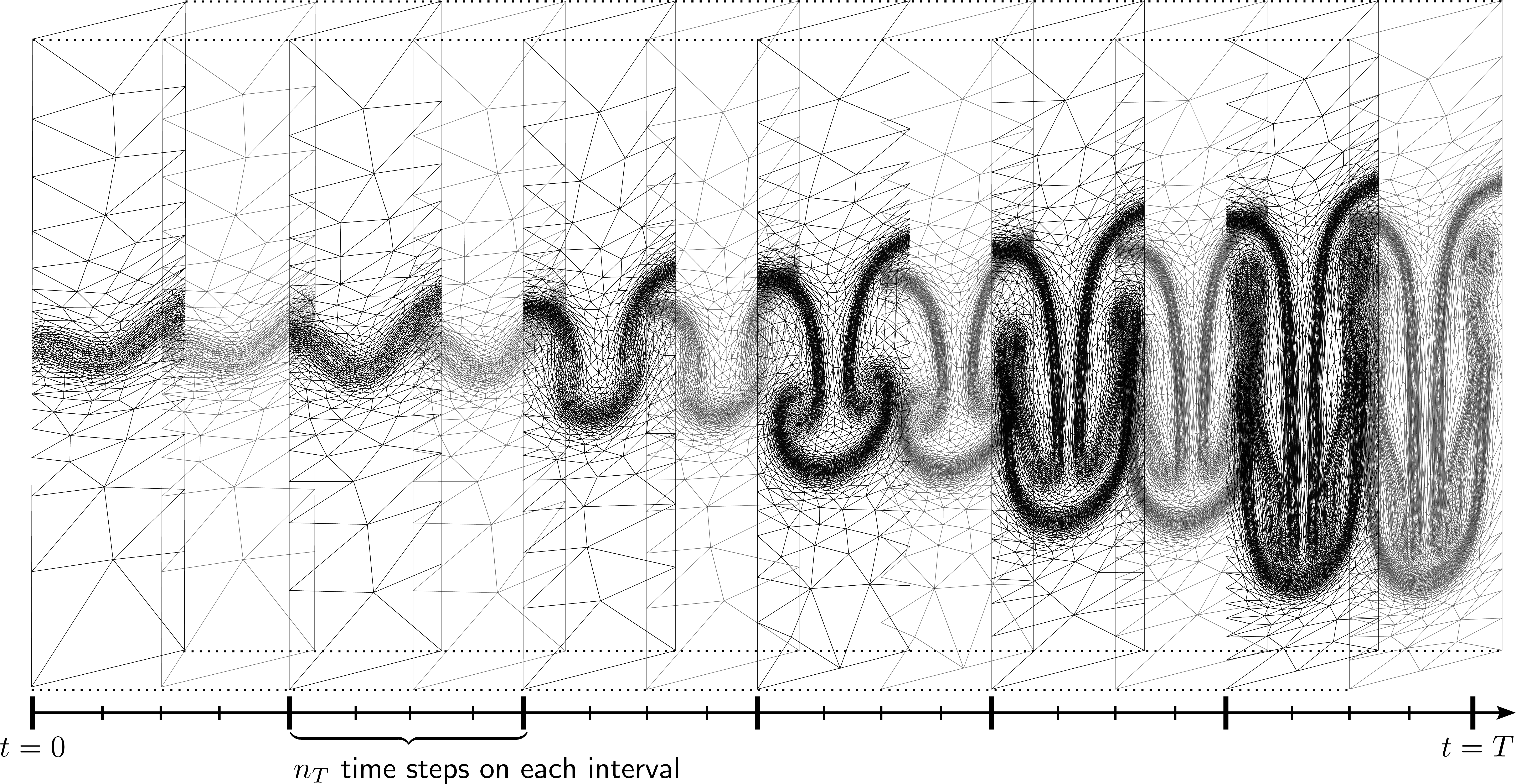}
  \caption{Space-time mesh adapted to a Rayleigh-Taylor instability.
  The time axis is sliced into $\nI$ sub-intervals, on which $\nt$ time steps each are computed.}
  \label{fig:st_mesh}
\end{figure}

For unsteady problems, we consider the space-time computational domain $\Omega_t \triangleq \Omega \times [0,T]$,
with $T > 0$ the (finite) maximum simulation time.
A mesh (or a collection of meshes) of $\Omega$,
together with a discretization of the time axis,
can be seen as a $d+1$-dimensional space-time mesh of $\Omega_t$,
and metric-based mesh adaptation on $\Omega_t$ amounts to determining a field of space- and time-varying
metrics $\mcal(\bfx,t)$ describing the optimal anisotropic mesh at each instant,
then generating a unit space-time mesh for this prescription.
In \citep{alauzet2018time}, the metric field minimizing the spatial error estimate integrated in time, given
hereafter in \Cref{sec:global_fp}, is derived,
assuming
a remeshing on equally spaced time sub-intervals.
Thus, the time-varying metric field is approximated by a collection of metrics $\mcal_i(\bfx)$,
which each summarize the behavior of the sensor (here, $\phi$) during the associated time sub-interval.
An example of $2+1$-dimensional space-time mesh adapted for a
Rayleigh-Taylor instability
is shown of \Cref{fig:st_mesh}.
In this paper, time discretization is done at constant time step, which translates
in equally-spaced slices in the depiction of the space-time mesh.

For space-time meshes, the \emph{space-time complexity} $\ccal_{st}$ is defined by
\begin{equation}
  \ccal_{st}(\mcal) \triangleq \int_0^T (\Delta t)^{-1}\left( \intom{\sqrt{\det \mcal(\bfx,t)}}\right)\,dt,
  \label{eq:st_complexity}
\end{equation}
where $\Delta t$ is the time step at time $t$.
Similarly to the spatial complexity, $\ccal_{st}$ is an estimate of the total number
of \emph{space-time} vertices,
that is, an estimate of the total resources used for an unsteady simulation.

\subsection{Global fixed-point mesh adaptation method}
\label{sec:global_fp}
For the unsteady adaptive procedure, we implemented the global
fixed-point method proposed in \citep{alauzet2018time}, based on a time-varying Riemannian metric.
The idea is to converge the mesh-solution pair together, accounting for the unsteadiness of the solution.
To this end, a first complete simulation is performed over the time interval $[0,T]$.
Then, knowing the behavior of the solution, the mesh is refined at equally
spaced time locations, and a more accurate solution is recomputed on the adapted meshes.
This process is repeated until convergence of the mesh-solution pair, typically a few (2-3) iterations.
More precisely, the time interval is partitioned into $\nI$ sub-intervals:
\begin{equation}
  [0,T] = [0 = t_1, t_2] ~\cup~ \ldots ~\cup~ [t_i, t_{i+1}] ~\cup~ \ldots ~\cup~ [t_{\nI},t_{\nI+1} = T].
  \label{eq:spatialErrorModel}
\end{equation}
On each sub-interval, the mesh is kept constant and $\nt$ solver time steps are
performed to advance the CHNS solution from $t_i$ to $t_{i+1}$.
Thus, $\nI$ metric fields $\mcal_i(\bfx)$ and adapted meshes $\tcal_{i}$ are generated for each global fixed-point iteration.
The complete fixed-point method is described by Algorithm 1 in \citep{alauzet2018time}. 
In this work, implicit time-stepping is used and time integration is done at constant time step $\Delta t = T/(\nI\nt)$.


The sensor driving the adaptation process is
a spatial interpolation error model for $\phi$ evaluated in $L^1([0,T]), L^p(\Omega)$ norm.
Assuming a linear representation for $\phi_h$, it writes in terms of the metric:
\begin{equation}
  \ecal^p(\mcal(\bfx,t)) = \int_0^T \int_\Omega \trace \left( \mcal^{-\frac{1}{2}}(\bfx,t) |H_\phi(\bfx,t)| \mcal^{-\frac{1}{2}}(\bfx,t)\right)^p\,d\bfx\,dt,
\end{equation}
where the matrix absolute value $|H_\phi|$ and exponent $\mcal^{-\frac{1}{2}}$ are applied to the eigenvalues, that is, $|A| \triangleq P|\Lambda|P^T$
and $A^{-\frac{1}{2}} \triangleq P\Lambda^{-\frac{1}{2}}P^T$ respectively.
The integrand is a local spatial error model based on a Taylor expansion, which involves derivatives of order $k+1$ when interpolating with degree $k$ polynomials.
In particular, this Hessian-based model is valid only for a linear approximation of $\phi$, however, higher-order models also exist \cite{coulaud2016very, bawin2024metric}.
The optimal Riemannian metric $\mcal_i$ on each time sub-interval is derived by minimizing
$\ecal^p(\mcal(\bfx,t))$ under the constraint $\ccal_{st}(\mcal) = \ncalst$ on the space-time complexity,
where $\ncalst$ is prescribed by the user.
This constraint controls the number of vertices in each adapted mesh, and thus reflects the resources allocated to the unsteady simulation.

\paragraph{Remark} It is important to note that $\ecal^p(\mcal)$ is a spatial error, albeit integrated in time.
Thus, the time discretization error is not \emph{a priori} controlled.
When an explicit time-stepping method is used, the spatial error also controls
the time discretization error through a CFL condition (see \citep{alauzet2018time} and the references therein).
Here, an implicit scheme is used and both errors are decoupled,
so we instead choose a small enough time step (through the product $\nI\nt$) to
yield a sufficiently small time discretization error.
A fully automatic approach with both time and space error control tailored for
implicit time-stepping scheme was proposed in \citep{sauvage2024space}, and will be
investigated for unsteady CHNS simulations in an upcoming work.



\subsubsection{Metric on each time sub-interval}
\label{sec:unsteadyMetric}
The problem of finding the set of metrics $\mcal_i(\bfx)$ minimizing \eqref{eq:spatialErrorModel} under the constraint $\ccal_{st}(\mcal) = \ncalst$ was tackled in \citep{alauzet2018time}.
On the $i$-th time sub-interval, the optimal metric is the scaled time integral of the Hessian on this sub-interval, and writes:
\begin{align}
  &\mcal_i(\bfx) = \ncalst^{\frac{2}{d}}
  \left( \sum_{j = 1}^{n_I} \mathcal{K}_j \left( \int_{t_j}^{t_{j+1}}(\Delta t)^{-1} \,dt \right)^{\frac{2p}{2p+d}}\right)^{-\frac{2}{d}}\\
  &~~~~~~~~~~~~~~~~~~~~\times
  \left( \int_{t_i}^{t_{i+1}}(\Delta t)^{-1} \,dt\right)^{-\frac{2}{2p+d}}
  (\det H_{\phi,i}(\bfx))^{-\frac{1}{2p+d}}
  H_{\phi,i}(\bfx),
\end{align}
where $d$ is the space dimension, $p$ refers to the minimized $L^p$ norm, and:
\begin{equation}
  H_{\phi,i}(\bfx) \triangleq \int_{t_i}^{t_{i+1}} |H_\phi(\bfx,t)| \,dt,
  ~~~~~\mathcal{K}_j \triangleq \intom{\left( \det H_{\phi,j}(\bfx) \right)^{\frac{p}{2p+d}}}.
\end{equation}
$H_{\phi,i}(\bfx)$ captures the evolution of the Hessian over the $i$th time sub-interval,
whereas the real number $\kcal_j$ is a global scaling factor on time interval $j$.
Following \citep{alauzet2018time}, the prescribed space-time complexity is $\ncalst = \ncalavg\nI,$
where $\ncalavg$ is the target average spatial complexity on each time interval.
Moreover, the time step is constant and given by $\Delta t = T/(\nI\nt)$, and all time intervals have the same length $t_{i+1}-t_i = T/\nI$, thus
$
  \int_{t_i}^{t_{i+1}}(\Delta t)^{-1} \,dt = \nt,
$
yielding the simplified metrics:
\begin{equation}
  \mcal_i(\bfx) =
  \left(\frac{\ncalavg\nI}{\nt}\right)^{\frac{2}{d}} \left(\sum_{j = 1}^\nI \mathcal{K}_j\right)^{-\frac{2}{d}} (\det H_{\phi,i})^{-\frac{1}{2p+2}} \, H_{\phi,i},
  \label{eq:metrics}
\end{equation}
which is the collection of metrics considered in the rest of the paper.

\subsection{Derivatives recovery}
The Hessian $H_\phi(\bfx,t)$ must be evaluated numerically to compute the
metrics \eqref{eq:metrics} on each interval. Various recovery methods exist to this end,
and we use Zhang and Naga's Polynomial Preserving Recovery (PPR) \cite{zhang2005new}
to reconstruct a quadratic approximation of $\phi$ by local least squares on
patches of elements, yielding a linear gradient.
The procedure is applied a second time to each component of the gradient to recover a linear Hessian.
To mitigate the effect of anisotropy of the recovered fields and obtain a better
conditioning of the least-squares system, a scaling is applied on each patch of triangles surrounding a vertex.
The complete procedure is detailed in \citep{briffard2017contributions, bawin2024metric}.

\subsection{Metric gradation}
Riemannian metrics for anisotropic mesh adaptation are typically derived from
fields exhibiting steep gradients or discontinuities,
and from numerical data which can be noisy. Here, the considered
phase marker field $\phi$ is continuous, but varies rapidly at the fluid-fluid interface,
and the Hessian is recovered by through two least-squares fittings.
This typically results in metric fields with abrupt variations in the prescribed sizes
and/or orientations, which can lead to the anisotropic mesh generator being unable
to produce a triangulation complying with the metric.
In practice, it is necessary to apply a gradation to the metric fields, to smooth
them and limit their variations.
Here, we apply the gradation method described in \cite{alauzet2010size}:
given a parameter $\beta > 1$,
the gradation procedure ensures that the length ratio of two consecutive mesh edges
is bounded by $\beta$.
This is enforced by spanning a metric field $\mcal_{\bfp}$ with prescribed size growth $\beta$
from each mesh vertex $\bfp$, then intersecting this metric with the size prescription
from all other vertices.
The iterative edge-based algorithm proposed in \cite{alauzet2010size} is linear
in the number of mesh edges and guarantees a smooth size variation bounded by $\beta$ throughout the domain.
Moreover, the edges are treated in a random order to limit the impact of the mesh
in the resulting size field.
In the numerical applications of Section \ref{sec:numericalApplications},
we typically set $\beta \in [1.5,2]$.

\subsection{Solution transfer from one interval to the next}
Lastly, as the domain is completely remeshed from one time sub-interval to the next,
the solution needs to be transferred from the mesh $\tcal_i$ to the next mesh $\tcal_{i+1}$.
In this work, this is done with a simple interpolation of the same order as the
field to be transferred (quadratic for $\bfu_h$, linear for $p_h, \phi_h, \mu_h$).
This transfer scheme is simple but not conservative,
hence using a large number of time intervals $\nI$ will eventually degrade the
quality of the numerical solution by accumulation of transfer error,
as discussed in \Cref{sec:mms_intervals} hereafter,
and conservative transfer schemes will be considered in an upcoming work.

\section{Numerical applications}
\label{sec:numericalApplications}

In this section, we verify the implementation of both the Cahn-Hilliard Navier-Stokes
solver and of the mesh adaptive procedure, then we apply them to the well-known rising bubble benchmark from Hysing et al. \citep{hysing2009quantitative}.
The simulations are performed with an in-house finite element solver\footnote{\url{https://github.com/arthurbawin/feNG}},
which takes care of the finite element resolution, derivatives recovery and computation and smoothing
of the metric fields.
Remeshing is performed with the open-source mesh generator
\texttt{mmg2d}\footnote{\url{https://github.com/MmgTools/mmg}} \citep{dobrzynski2012mmg3d}.
The linear system at each Newton iteration is solved with Intel's MKL PARDISO multithreaded direct solver.


\subsection{Verification of the CHNS solver}
The accuracy of the solver is assessed using two sets of manufactured solutions.
For all convergence tests, the simulation time is $[0,T]$ with $T = 1$
and the model parameters are set to $\rho_1 = 1000, \rho_2 = 1,
\visc_1 = 100, \visc_2 = 1, \sigma = 0.01, \epsilon = 10^{-3}$ and $M = 10^{-6}$.
The error is computed using an $L^1$ norm in time and the relevant spatial norm:
\begin{align}
  E = \Vert e \Vert_{L^1([0,T]), V(\Omega)}
  &= \int_0^T \Vert e(\bfx,t) \Vert_{V} \,dt \simeq \Delta t \sum_{i = 1}^{N_t}  \Vert e(\bfx,t_i) \Vert_{V},
\end{align}
where $V = (H^1(\Omega))^2$ for the velocity, $L^2(\Omega)$ for the pressure, and $H^1(\Omega)$ for the phase marker and the chemical potential.

\subsubsection{Space convergence}
To assess the spatial accuracy, a convergence study is performed
on a sequence of uniformly refined triangulations described by their characteristic mesh size $h_n$ for $n = 0 \to 7$.
We select the following fields:
\begin{equation}
\begin{alignedat}{1}
  \bfu(\bfx,t) &= t^4 \left(\sin(\pi x)\sin(\pi y), \cos(\pi x)\cos(\pi y) \right),\\
  p(\bfx,t) &= \sin(\pi x)\cos(\pi y),\\
  \phi(\bfx,t) &= t^4 \sin(\pi x)\sin(\pi y), \, ~~~~~~\mu(\bfx,t) = \sin(\pi x)\sin(\pi y),
\end{alignedat}
\end{equation}
and add the corresponding source terms to
\cref{eq:CHNS_u,,eq:CHNS_phi,eq:CHNS_mu}.
The initial and boundary conditions are set with the exact solution,
and so is the first time step of the BDF2 scheme.
To isolate the effect of spatial discretization, we set the number of time steps $N_{t,n} = 100 \times 2^{n}$,
which is sufficiently high to have the spatial error dominate the time discretization error for all meshes.

\begin{figure}[t]
  \centering
  \includegraphics[width=\linewidth]{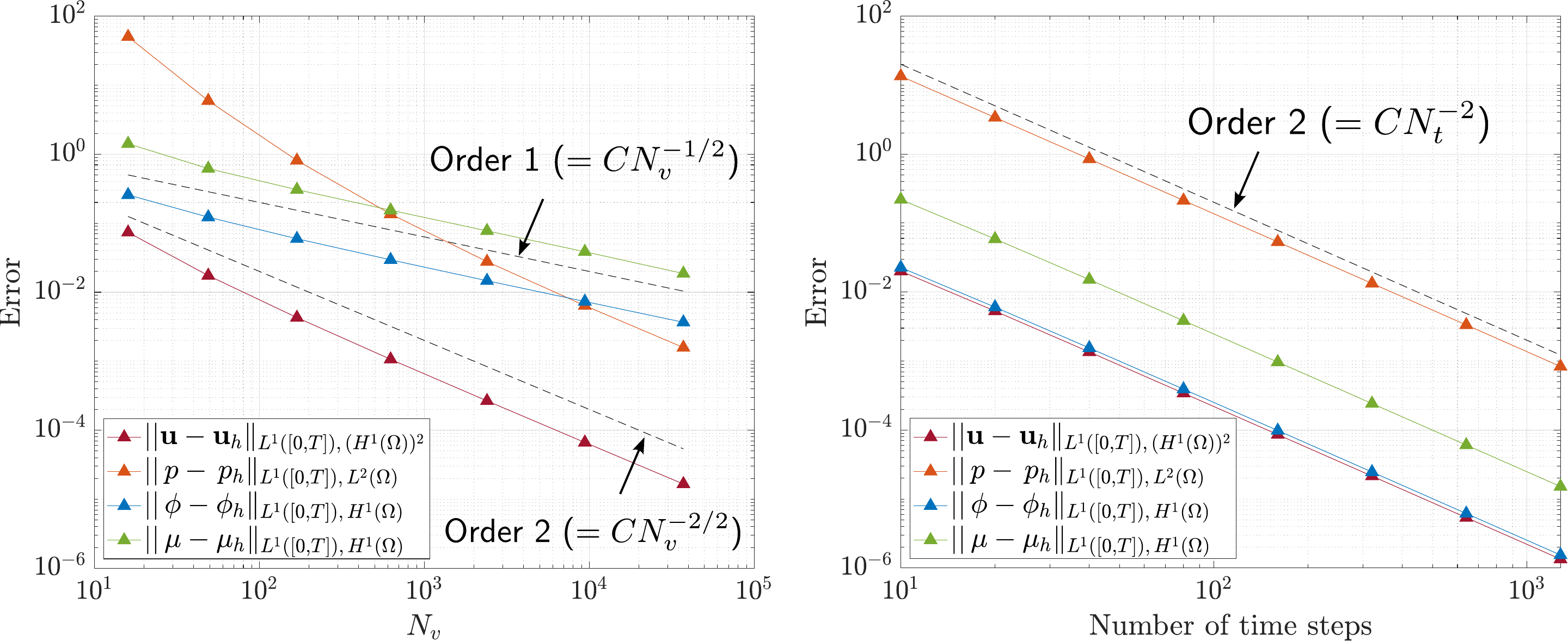}
  \caption{Verification of the CHNS solver: spatial (left) and temporal (right) convergence.}
  \label{fig:verificationCHNS}
\end{figure}

With Taylor-Hood elements,
we can expect
$\ocal(h^2)$ convergence for both the velocity error in $(H^1(\Omega))^2$ norm and the pressure error in $L^2$ norm.
Similarly, we expect $\ocal(h)$ convergence for the phase marker and chemical
potential errors in $H^1(\Omega)$ norm with linear elements.
The optimal spatial rates are obtained for $\bfu, \phi$ and $\mu$, left plot of \Cref{fig:verificationCHNS}.
The convergence rate for the pressure is higher than expected, but decreases towards 2 in the asymptotic regime.

\paragraph{Remark}
When working with anisotropic meshes,
a unique size $h$ no longer characterizes the fineness of the mesh,
as the element size typically varies greatly throughout the domain.
Instead, convergence studies are typically written in terms of either
the number of mesh vertices $N_v$ or mesh elements $N_e$.
On \emph{isotropic} meshes, the asymptotic relation $N_v \sim h^{-d}$ holds between
the mesh size and the number of mesh vertices, so that convergence results can be
reported as the error in terms of $h \sim N_v^{-1/d}$.
Convergence with a rate $r$ is thus
characterized by a slope $-r/d$ in a $(\log E, \log N_v)$ graph.
For consistency, convergence results on both isotropic (such as here) and
anisotropic meshes (in \Cref{sec:mms_adaptive}) are reported in terms of the number of mesh vertices $N_v$
 (or, for \Cref{sec:mms_adaptive}, the number of mesh vertices in the space-time mesh, see hereafter).

\subsubsection{Time convergence}
The accuracy of the implicit time integrator (BDF2) is verified by refining the time step.
To isolate the effect of the time discretization error, we select a manufactured solution
whose spatial component is interpolated exactly with Taylor-Hood and equal order linear elements,
thus a quadratic velocity and linear $p, \phi$ and $\mu$:
\begin{equation}
\begin{alignedat}{1}
  \bfu(\bfx,t) &= t^4 \left(x^2 - 2x + 3y^2 + y - xy + 1, 3x^2 + y^2/2 + 2y - 1 - 2xy\right)\\
  p(\bfx,t) &= 10g(x,y), ~~~~~\phi(\bfx,t) = t^4 g(x,y), ~~~~~\,\mu(\bfx,t) = g(x,y),
\end{alignedat}
\end{equation}
with $g(x,y) = (x + y - 1)$.
The expected second-order convergence is achieved, right plot of \Cref{fig:verificationCHNS},
here given in terms of the number of time steps $N_t$.

\subsection{Verification of the metric field for unsteady anisotropic adaptation}
\label{sec:mms_adaptive}
Next, we perform a space-time convergence study to assess the implementation of
the global transient fixed-point method.
The goal of this verification is to show that the unsteady adaptation process
exhibits at least second-order convergence, as indicated in \citep{alauzet2018time},
when adapting in space with the optimal space-time metric and using a implicit
time integration scheme with constant time stepping.
We consider on $\Omega \times [0,T] = [-2,2] \times [-1,1] \times [0,T]$ the manufactured solution:
\begin{equation}
  \phi_\text{ref}(\bfx,t) = \tanh\left( \frac{2(x -ct) - \sin 5y}{\delta} \right),
\end{equation}
with $c = 1$ and $\delta = 0.02$.
This reference solution is an oscillating hyperbolic tangent traveling across $\Omega$ at unit velocity $c$, \Cref{fig:tanh}, and whose steepness is controlled by $\delta$.
The solved PDE is the unsteady heat equation with initial and boundary conditions:
\begin{equation}
  \label{eq:MMS_heat}
  \begin{alignedat}{2}
  \ddt{\phi} - \Delta\phi + f &= 0 &&\text{ on }\Omega,\\
  \phi(\bfx,0) &= \phi_\text{ref}(\bfx, 0) ~~&&\text{ on }\Omega,\\
  \phi(\bfx,t) &= \phi_\text{ref}(\bfx, t) &&\text{ on }\partial\Omega,
  \end{alignedat}
\end{equation}
with $f$ set so that $\phi_\text{ref}$ is solution of the system.
The discrete solution $\phi_h$ is computed with linear elements,
and we consider the interpolation error $e = \phi_\text{ref} - \phi_h$ in $L^1([0,T]), L^2(\Omega)$ norm:
\begin{equation}
  E = \Vert e \Vert_{L^1([0,T]), L^2(\Omega)}
  \simeq \sum_{j = 1}^{\nI}\nt \Delta t \, \Vert e(\bfx,t_{j})\Vert_{L^2(\Omega)}.
\end{equation}
As in \citep{alauzet2018time}, we examine the convergence of the error $E$ with respect to the \emph{effective} space-time complexity $\nst$, defined by
\begin{equation}
  \nst = \sum_{i = 1}^\nI \nt \times N_{v,i} = \nt \sum_{i = 1}^\nI N_{v,i},
\end{equation}
where $N_{v,i}$ is the number of vertices on the $i$-th adapted mesh.
Since $\nt$ time steps are computed on $\nI$ meshes of $N_{v,i}$ vertices each, $\nst$ reflects the total allocation of resources for the simulation.

In the following, the overall error $E$ is plotted against $\nst$,
and we discuss the influence of the time step,
the prescribed average mesh density $\ncalavg$
and the number of time sub-intervals $\nI$.
Given that we consider the complexity of a $d+1$-dimensional \emph{space-time} mesh,
convergence is computed in \citep{alauzet2018time} as $E \sim \nst^{-r/(d+1)} = \nst^{-r/3}$, with $r$ the convergence rate.
This assumes that all $d+1$ dimensions are refined, which is the case in \citep{alauzet2018time} since an explicit time-stepping scheme is used with CFL control.
Here, the time step is fixed throughout a convergence study, so that only the spatial part of the mesh is refined.
Consequently, in the following, convergence is presented in terms of $\nst^{-r/d}$ with $d = 2$.
To generate each data point, \eqref{eq:MMS_heat} is solved until $T = 1$ using $n_I$ time sub-intervals,
and 5 fixed-point iterations are performed, that is, the simulation is run 5 times, starting on a coarse uniform mesh.
All tests are performed at constant time step $\Delta t = T/(\nI\nt)$,
which is achieved by either keeping both $\nI$ and $\nt$ fixed,
or decreasing $\nt$ when $\nI$ increases.
Lastly, the time integral of the Hessian $H_{\phi,i}$ is computed by evaluating the absolute Hessian matrix at each time step.

\paragraph{Remark}
Here, the interface represented by the manufactured solution is simply advected
throughout $\Omega$ and does not undergo any particular deformation.
As a result, each adapted mesh $\tcal_i$ contains roughly the same number $N_{v,i} = C\ncalavg/\nt$ of vertices,
with $C$ typically around 1.1. Thus, the effective complexity is
\begin{equation}
  \nst = \nt \sum_{i = 1}^\nI N_{v,i} = C\nI\ncalavg = C\ncalst
\end{equation}
and is also roughly the same as the prescribed complexity.
When the interface undergoes deformation, the number of vertices for each mesh varies,
but we still have $\sum_i N_{v,i} = C\nI\ncalavg/\nt.$

\paragraph{Remark}
From the expression of the metric \eqref{eq:metrics}, it is worth noting that
increasing $\nt$ at fixed $\ncalavg$ yields a lower mesh
density, thus a larger spatial error.
In terms of space-time complexity, advancing $\nt$ time steps on a mesh with $\ncalavg/\nt$
always yields a complexity of $\ncalavg$.
However, while complexity is multiplicative, errors are additive, that is,
setting a small time steps
on a poorly resolved mesh yields a large spatial error, and vice versa,
even though these configurations have comparable complexities.
Thus, for convergence studies, it is not beneficial to increase
the number of time steps $\nt$ at $\ncalavg$ fixed, as the spatial error rapidly
takes over and spoils the convergence rates.

\subsubsection{Influence of the time step and average mesh density at fixed $n_I$}
\label{sec:mms_fixed_ni}
We first present convergence results when the average mesh density $\ncalavg$ is increased
and for several values of $\Delta t$,
keeping constant both the number of time sub-intervals $n_I$ and the number of time steps per interval $n_T$.
For each convergence study, the prescribed average mesh density varies from $\ncalavg = 1,000$ to $128,000$,
for different combinations of $n_I = [4,16]$ intervals and $n_T = [2, 5, 10, 20]$ time steps per interval.
The prescribed space-time complexity thus ranges from $\ncalst = 4,000$ to $2.05$ million,
and the time step varies from $\Delta t = 1/8$ to $1/160$.
In practice, due to metric gradation and the specificities of the mesh generator,
we generally observe about 10\% more vertices than prescribed,
and the adapted meshes contain between 1,100 and 138,000 vertices.

\begin{figure}[t]
  \centering
  \includegraphics[width=0.65\linewidth]{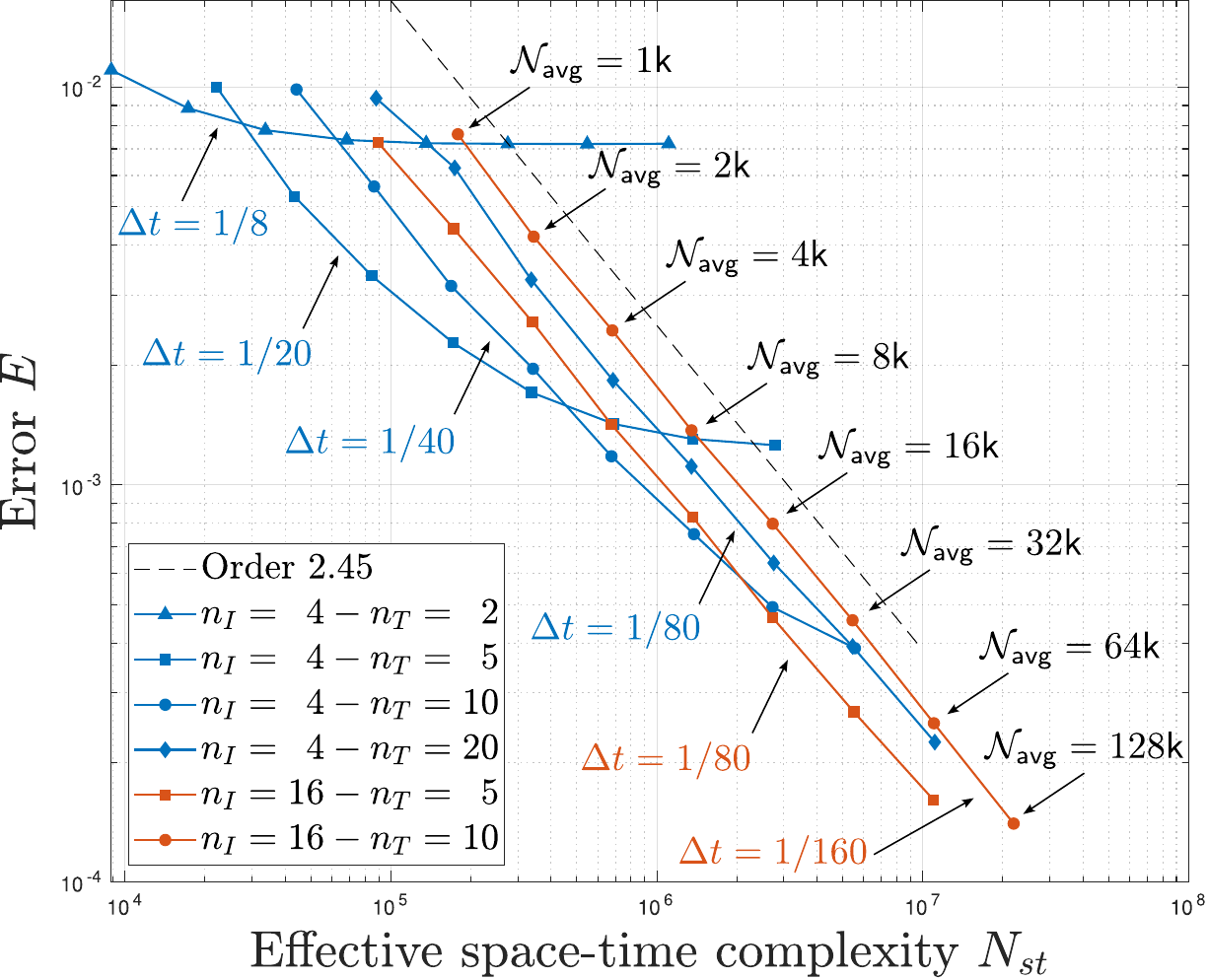}
  \caption{Convergence study for the manufactured problem \eqref{eq:MMS_heat}.
  The time step decreases from one convergence curve to the next, but is kept constant on a given convergence curve.
  Evolution of the error for $\nI = 4$ or $16$ intervals as the average mesh density $\ncalavg$ increases.}
  \label{fig:verificationST_fixedIntervals}
\end{figure}

The evolution of $E$ vs $\nst$ is given in \Cref{fig:verificationST_fixedIntervals}.
For larger time steps $(\geq 1/40)$,
the overall error decreases until it becomes controlled by its time component and reaches a plateau.
A convergence at order 1.6 is obtained with small enough time steps (the black dashed line is $C\nst^{-1.6/2}$).


\subsubsection{Influence of the number of time sub-intervals at fixed $\ncalavg$}
\label{sec:mms_intervals}
\begin{figure}[h!]
  \centering
  \includegraphics[width=\linewidth]{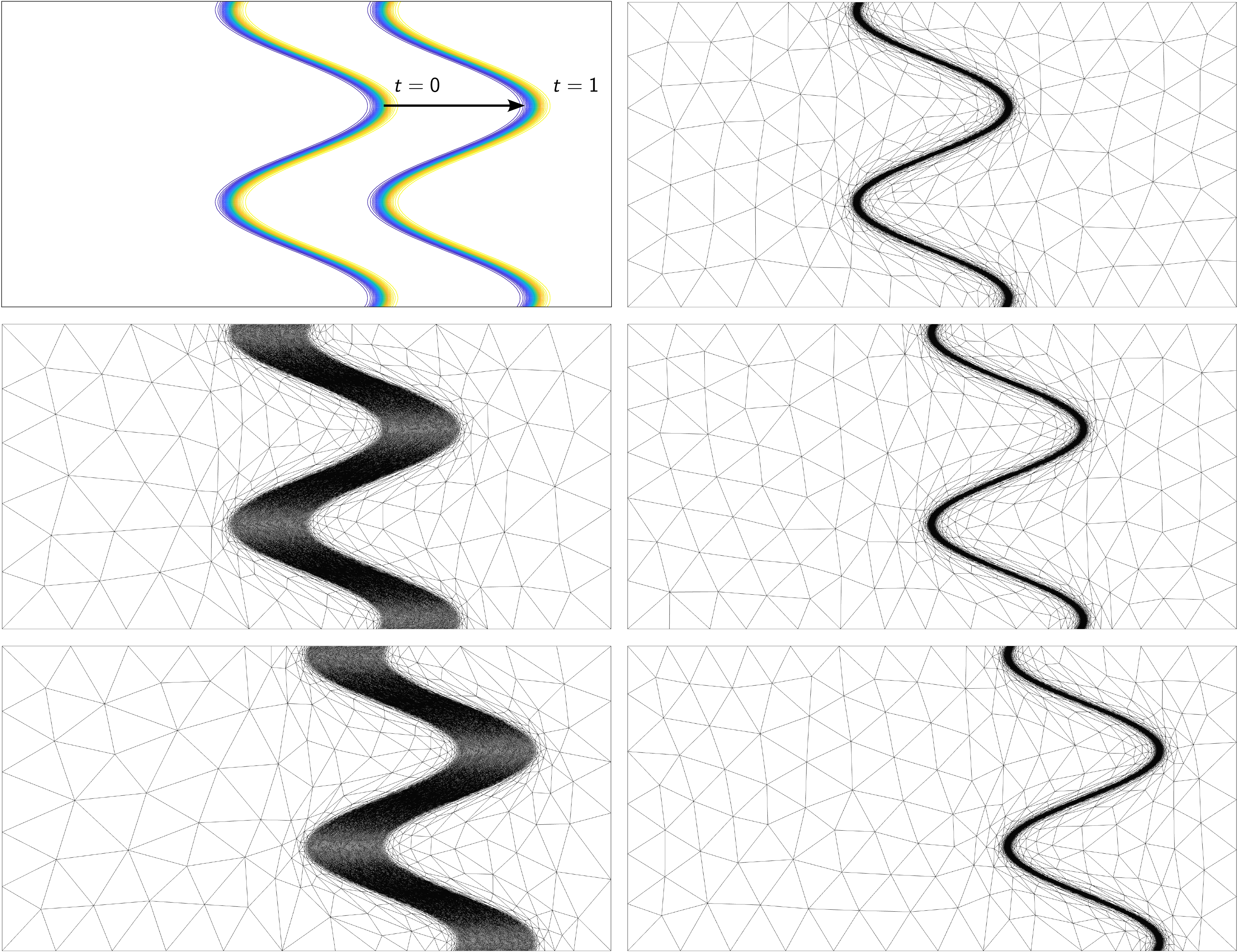}
  \caption{Top left: the manufactured solution $\phi_\text{ref}$ at $t = 0$ and $t = 1$.
  Middle and bottom left: adapted meshes for the time sub-intervals $1$ and $2$ (out of 2)
  for a prescribed mesh density $\ncalavg = 6.4$M.
  These meshes contain 13,965 and 13,987 vertices respectively.
  Right: adapted meshes for the intervals $1$, $64$ and $128$ (out of 128) and $\ncalavg = 100$k.
  These meshes contain 13,341, 13,254 and 13,321 vertices respectively.
  }
  \label{fig:tanh}
\end{figure}

Here, the number of time sub-intervals is increased while keeping
the prescribed average mesh density $\ncalavg$ constant.
The time step is set to $\Delta t = T/(\nI\nt) = T/1024$,
smaller than the ones in \Cref{sec:mms_fixed_ni} to ensure
that the overall error is controlled by its spatial component.
Intuitively, increasing the number of intervals
is beneficial as it limits the complexity required to reach a certain error level.
For instance, an adapted mesh for $\phi_\text{ref}$ with $\nI = 1$
will be refined along the whole part of the domain swept by the tanh.
But at each instant, only the mesh elements along the steep gradient are really
required, as refined elements in constant regions
add mesh complexity for no error reduction.
This is represented on \Cref{fig:tanh}: on the left are meshes with $\ncalavg = 6.4$M
and 2 time intervals ($\nt = 512$ time steps), and on the right are meshes with
$\ncalavg = 100$k and 128 intervals ($\nt = 8$).
The prescribed complexity $\ncalst = 12.8$M is identical for both configurations,
and the complexities effectively obtained are $\nst = 14,311,424$ and $13,633,296$ respectively.
However, the space-time error is $E = 3.02 \times 10^{-2}$ for 2 intervals and
$1.91 \times 10^{-3}$ for 128 intervals.
Indeed, each mesh on the right contains roughly the same number of vertices as the
left ones, but concentrates these vertices around the steep gradient.

%

It is clear, however, that a larger number of intervals comes with more solution transfers,
which will eventually spoil the quality of the numerical solution.
Here, to isolate the effect of $\nI$ on the total error,
the transfer error is cancelled by reinterpolating the exact solution
at the beginning of each time interval.
This effectively sets the temporal error to zero at the beginning of each interval,
so that the only contribution to the total error is the spatial component,
and is consistent with choosing a time step small enough to assume a small time discretization error.


\begin{figure}[t!]
  \centering
  \includegraphics[width=\linewidth]{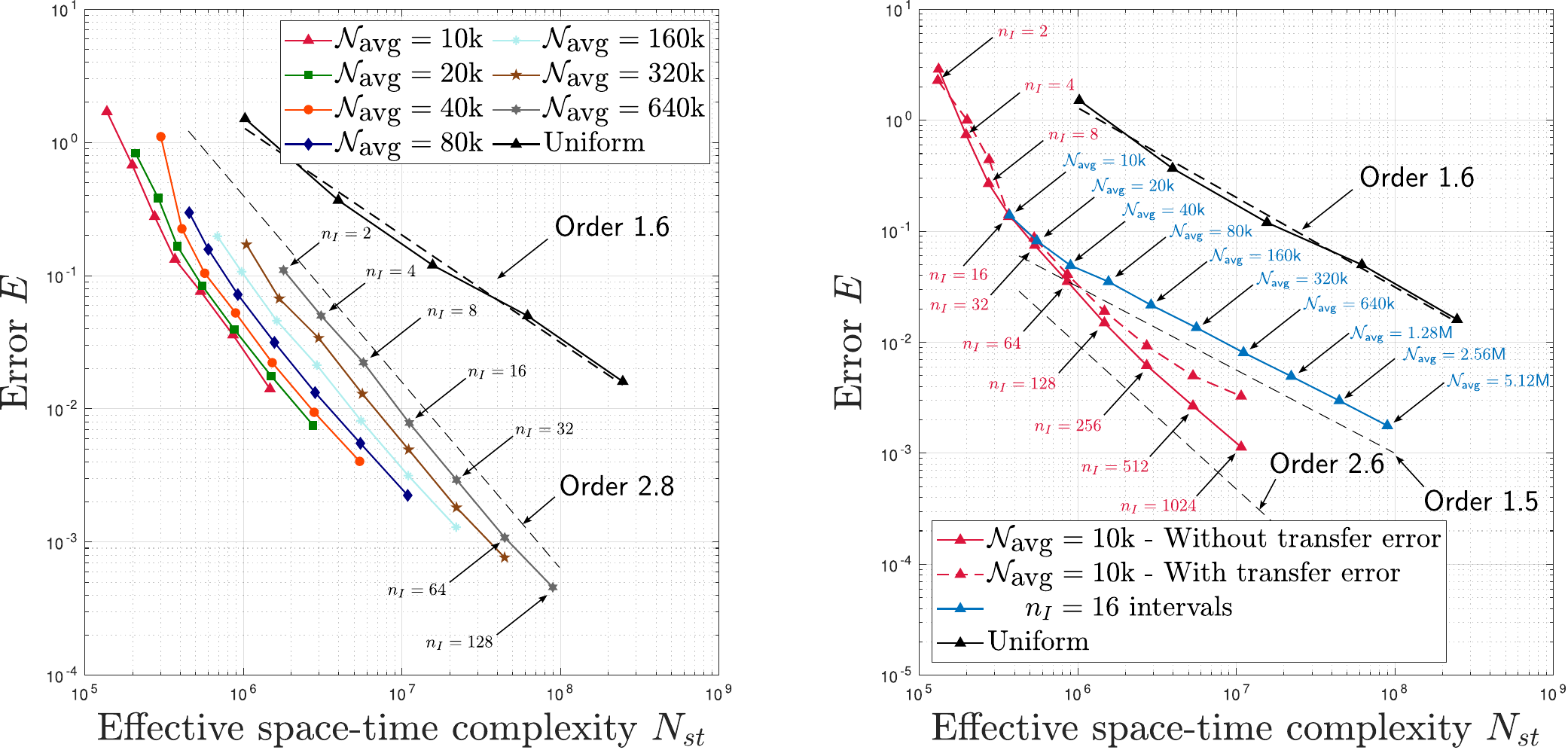}
  \caption{Convergence curves for the manufactured problem \eqref{eq:MMS_heat} at fixed time step $\Delta t = 1/1024$.
  Left: error for a fixed prescribed average mesh density $\ncalavg$ and increasing number of time sub-intervals.
  Right: comparison between error at fixed number of intervals and increasing average mesh density vs. fixed mesh density and increasing number of time intervals.}
  \label{fig:verificationST_fixedDensity}
\end{figure}

Convergence curves for an increasing mesh density $\ncalavg = 10$k to $640$k
are presented on the left plot of \Cref{fig:verificationST_fixedDensity}.
On each curve, the number of sub-intervals varies from $\nI = 2$ to $128$. 
A convergence at order 2.8 is observed ($\nst^{-2.8/2}$),
whereas convergence is only of order 1.6 with uniform meshes due to the presence of steep gradients.
For a given complexity, the data points for various $\ncalavg$ are mostly stacked on top of one another,
indicating that lower error is achieved for a larger number of intervals, rather than larger mesh density,
in agreement with the discussion above and the observations in \citep{alauzet2018time}.
This is emphasized on the right plot of \Cref{fig:verificationST_fixedDensity},
where the evolution of the interpolation error at fixed density is compared against a fixed number of time intervals.
Convergence at order 2.6 is obtained when the number of time intervals is increased,
and at order 1.5 when increasing mesh density.
This difference of convergence rate for constant $\nI$ is also reported in \citep{alauzet2018time},
albeit for discontinuous fields.
This is true as long as the number of intervals remains modest,
otherwise the transfer error accumulates as $\nI$ increases and greatly impacts
the convergence rate,
as shown by the difference between the plain and dashed red lines on the right plot.

\subsection{Rising bubbles}
We now apply the global fixed-point method to the two well-known rising bubbles
benchmarks proposed by Hysing et al. \citep{hysing2009quantitative}:
a circular bubble of lighter fluid is immersed in a column of heavier fluid
and rises under the effect of buoyancy. The shape of the bubble as it rises is
influenced by the density and viscosity ratios
of the fluids, the Reynolds number $Re = \rho_1\sqrt{gD_0^3}/\visc_1$, and the Eötvös number $Eo = \rho_1gD_0^2/\sigma$,
which quantifies the balance between gravity and surface tension forces.

\begin{figure}[h!]
  \centering
  \includegraphics[width=\linewidth]{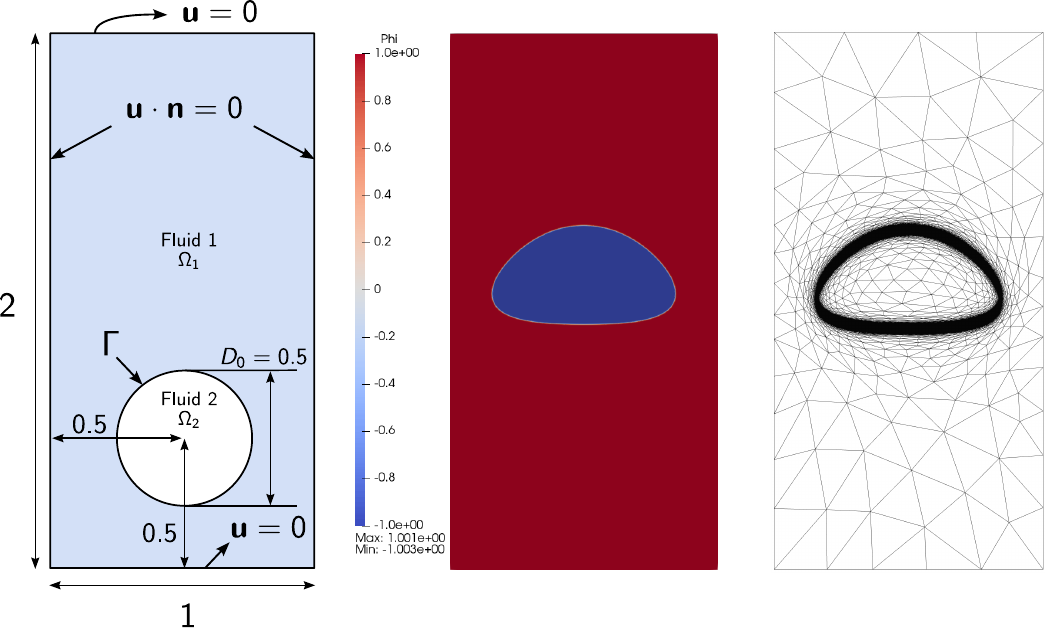}
  \caption{Left: Geometry of the rising bubble test cases \citep{hysing2009quantitative}.
  Middle and right: Phase field for benchmark 1 and $\epsilon = 0.00125$ at $t = 3$ s, and last adapted mesh for a target $\ncalavg = 4.1$M.
  The mesh contains 23,391 vertices and 46,746 triangles.}
  \label{fig:bubble_geometry_marker_mesh}
\end{figure}

The geometry is described on \Cref{fig:bubble_geometry_marker_mesh}: the bubble of fluid 2
and of radius $R_0 = 0.25$ is initially located at $(x_c,y_c) = (0.5,0.5)$, in a
column of fluid 1 of dimensions $[0,1] \times [0,2]$.
Slip boundary conditions are applied on the left and right walls,
and non-slip conditions are applied on the top and bottom.
The fluid parameters for both benchmarks are given in \Cref{tab:bubblesParameters}:
the first benchmark is characterized by a larger surface tension, so that
the bubble only undergoes mild deformation as it rises in the column.
Surface tension is much smaller for the second test case,
yielding a non-convex shape and filaments on both sides of the bubble.
These filaments may or may not separate from the rest of the bubble during the
simulation, depending on the adopted multiphase model
and discretization \citep{aland2012benchmark, hysing2009quantitative, brunk2025simple}.

The post-processed quantities available in the literature\footnote{\url{https://wwwold.mathematik.tu-dortmund.de/~featflow/en/benchmarks/cfdbenchmarking/bubble/bubble\_verification.html}}
for this benchmark are the bubble's area $A_b$, vertical center of mass $y_c$, vertical rising
velocity $u_b$ and circularity $C_b$.
These are respectively defined by:
\begin{equation}
  A_b= \int_{\Omega_2}\,d\bfx, ~~~y_c = \frac{1}{A_b}\int_{\Omega_2}y\,d\bfx, ~~~u_b = \frac{1}{A_b}\int_{\Omega_2}u_y\,d\bfx, ~~~ C_b = \frac{2\pi R_0}{\int_{\Gamma}\,d\ell},
\end{equation}
where the bubble domain is $\Omega_2 = \left\{ \bfx \in \Omega ~|~ \phi(\bfx) \leq 0\right\}$.

\begin{table}[t!]
\centering
\begin{tabular}{ccccccccc}
\toprule
Benchmark & $\rho_1$ & $\rho_2$ & $\visc_1$ & $\visc_2$ & $g$ & $\sigma$ & Re & Eo\\
\midrule
1 & 1000 & 100 & 10 & 1 & 0.98 & 24.5 & 35 & 10\\[0.6em]
2 & 1000 & 1   & 10 & 0.1 & 0.98 & 1.96 & 35 & 125\\
\bottomrule
\end{tabular}
\caption{Flow parameters for the rising bubble benchmarks \citep{hysing2009quantitative}.}
\label{tab:bubblesParameters}
\end{table}

\subsubsection{Adaptive strategy}
The simulation time for both benchmarks is $[0, T]$ with $T = 3$\,s,
and we set $\nI = 15$ time sub-intervals with $\nt = 50$ time steps each,
yielding a constant time step $\Delta t = 3/750 = 4 \times 10^{-3}$ s.
The global transient fixed-point strategy targets the minimization of
the $L^4$ norm of the interpolation error \eqref{eq:spatialErrorModel}, a value of $L^p$
known to give good results for most feature-based adaptation problems.

In addition to yielding an optimal mesh for a given interface thickness $\epsilon$,
anisotropic adaptation allows for a dynamic control of $\epsilon$.
For each simulation,
the initial value of $\epsilon$ is set to $\epsilon_0 = 0.64h$, similarly to \citep{brunk2025simple},
with $h$ the characteristic size of the starting uniform mesh.
The starting mesh for both benchmarks is a $32 \times 64$ uniform grid, yielding
$h = 1/32$ and $\epsilon_0 = 0.02$.
Two fixed-point iterations are performed to obtained a converged adapted mesh for this thickness,
then $\epsilon$ is divided by 2.
This process is repeated 4 times, performing 2 iterations to readapt the mesh each time,
yielding a final thickness of $\epsilon/16$ after a total of 10 fixed-point iterations.
This strategy allows to start from a relatively coarse uniform mesh and,
given a target complexity, obtain a  well-resolved interface in a completely automatic way.

\subsubsection{Results}
\begin{figure}[h!]
  \centering
  \includegraphics[width=\linewidth]{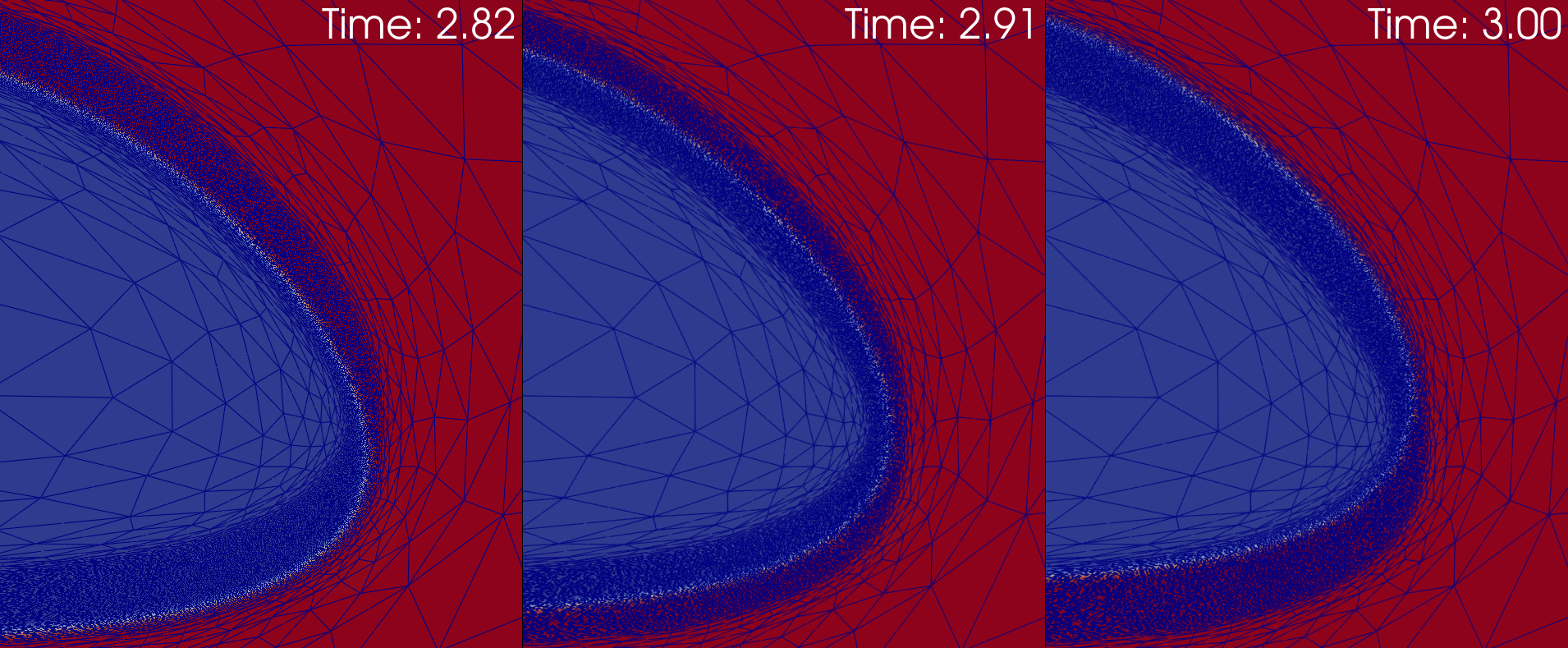}
  \caption{Bubble benchmark 1: Evolution of $\phi$ during the last time interval $[2.8,3]$.}
  \label{fig:bubble1_zoomMesh}
\end{figure}

\paragraph{Benchmark 1}
The phase marker field for the first benchmark with $\epsilon = \epsilon_0/16 = 0.00125$
at $t = 3\,$s is shown on \Cref{fig:bubble_geometry_marker_mesh},
together with the adapted mesh for the last time sub-interval $t \in [2.8, 3]$.
The adaptive process was run with a prescribed complexity $\ncalavg = 2.048$M,
and the last mesh contains 45,243 vertices and 90,457 triangles.
The over- and undershoot of the phase marker is limited to 0.3\%, as $\phi$ stays within $[-1.003, 1.001]$.
The mesh captures the movement of the bubble during the last time interval,
as shown on \Cref{fig:bubble1_zoomMesh}, which depicts the position of the interface
for three time steps of the last interval.

\begin{figure}[h!]
  \centering
  \includegraphics[width=\linewidth]{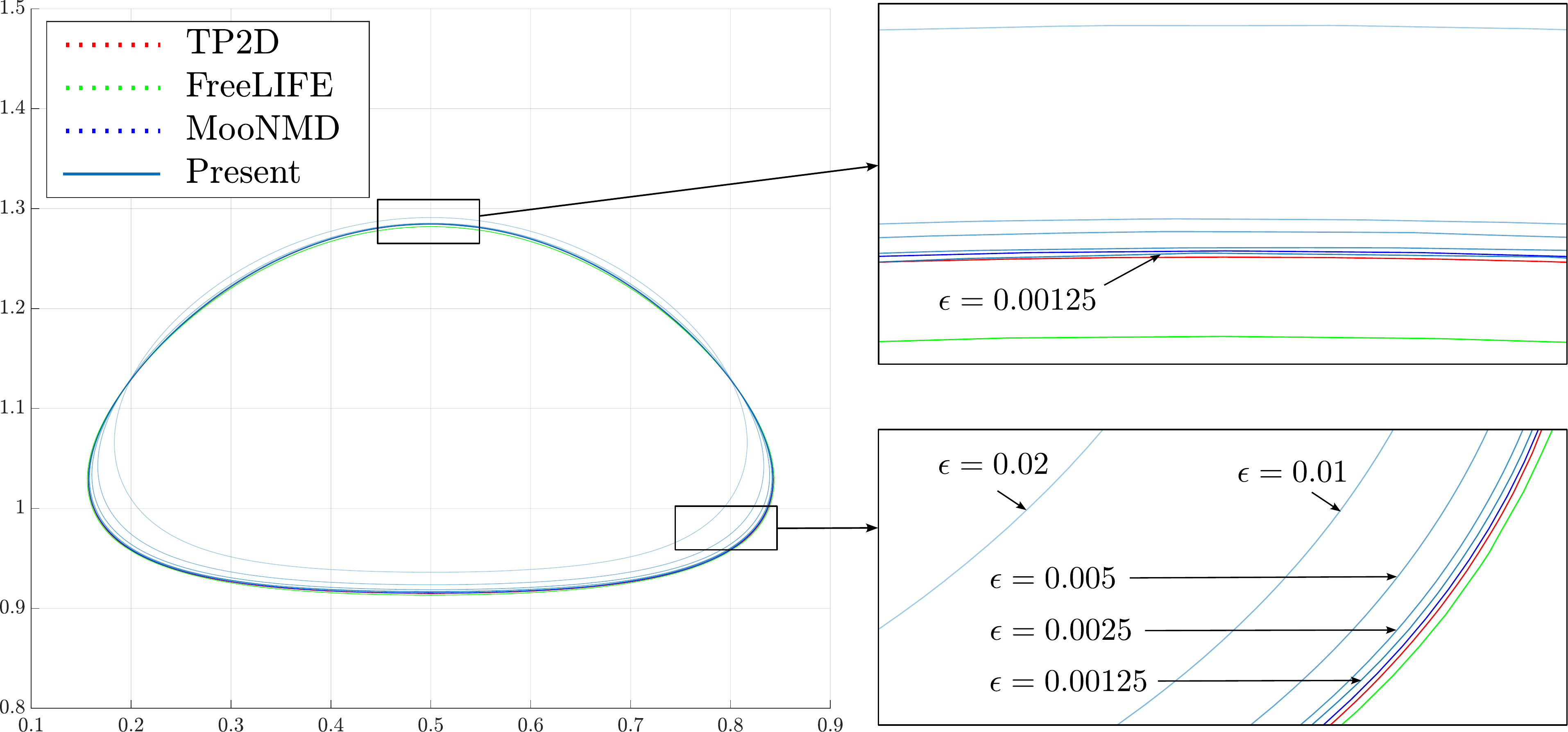}
  \caption{Bubble benchmark 1: Shape for decreasing interface thickness $\epsilon$.}
  \label{fig:bubble1_shape}
\end{figure}

The shape of the bubble for various interface thicknesses
is shown on \Cref{fig:bubble1_shape} (the opacity of the curves increases
as $\epsilon$ decreases) and is compared with
the results of TP2D, FreeLIFE and MooNMD on their finest discretizations.
The shapes for increasing $\epsilon$ converge towards the ones from the multiphase literature.
The agreement for $\epsilon = 0.00125$ is excellent, and the contour lies between the ones of TP2D and FreeLIFE.

\begin{figure}[hbtp!]
  \centering
  \includegraphics[width=\linewidth]{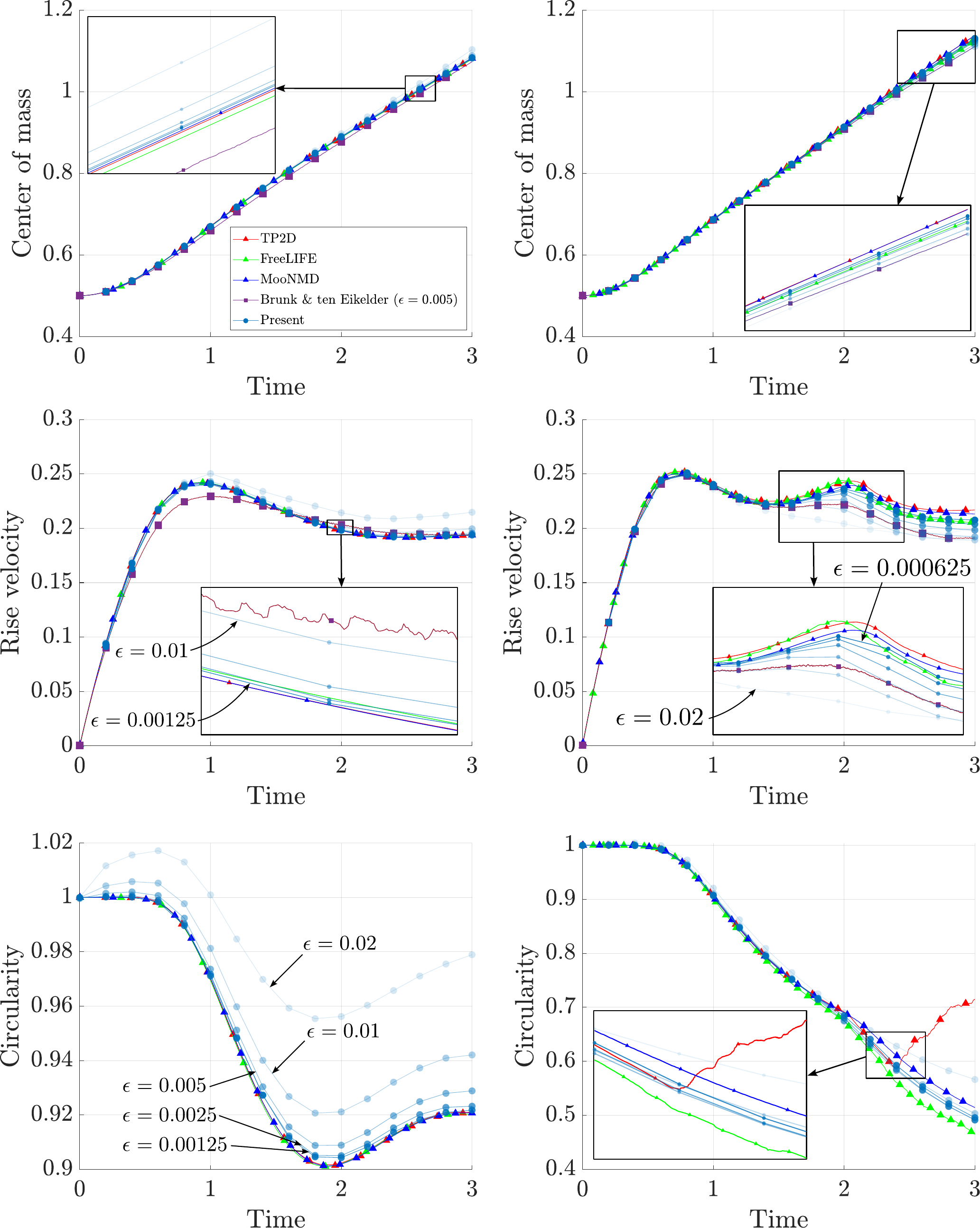}
  \caption{Bubble benchmark 1 (left) and 2 (right). Comparison with the literature: $y$-position of the center of mass, rise velocity and circularity.}
  \label{fig:bubblesIndicators}
\end{figure}

Similarly, we obtain a very good agreement for all the quantitative indicators,
left of \Cref{fig:bubblesIndicators}.
Our results are compared with those from
Hysing \citep{hysing2009quantitative},
as well as with CHNS simulation data at $\epsilon = 0.005$ from Brunk and ten Eikelder
\citep{brunk2025simple} ($y_c$ and $u_c$ only),
who used a CHNS formulation in terms of mass-averaged velocity.


\begin{figure}[h!]
  \centering
  \includegraphics[width=\linewidth]{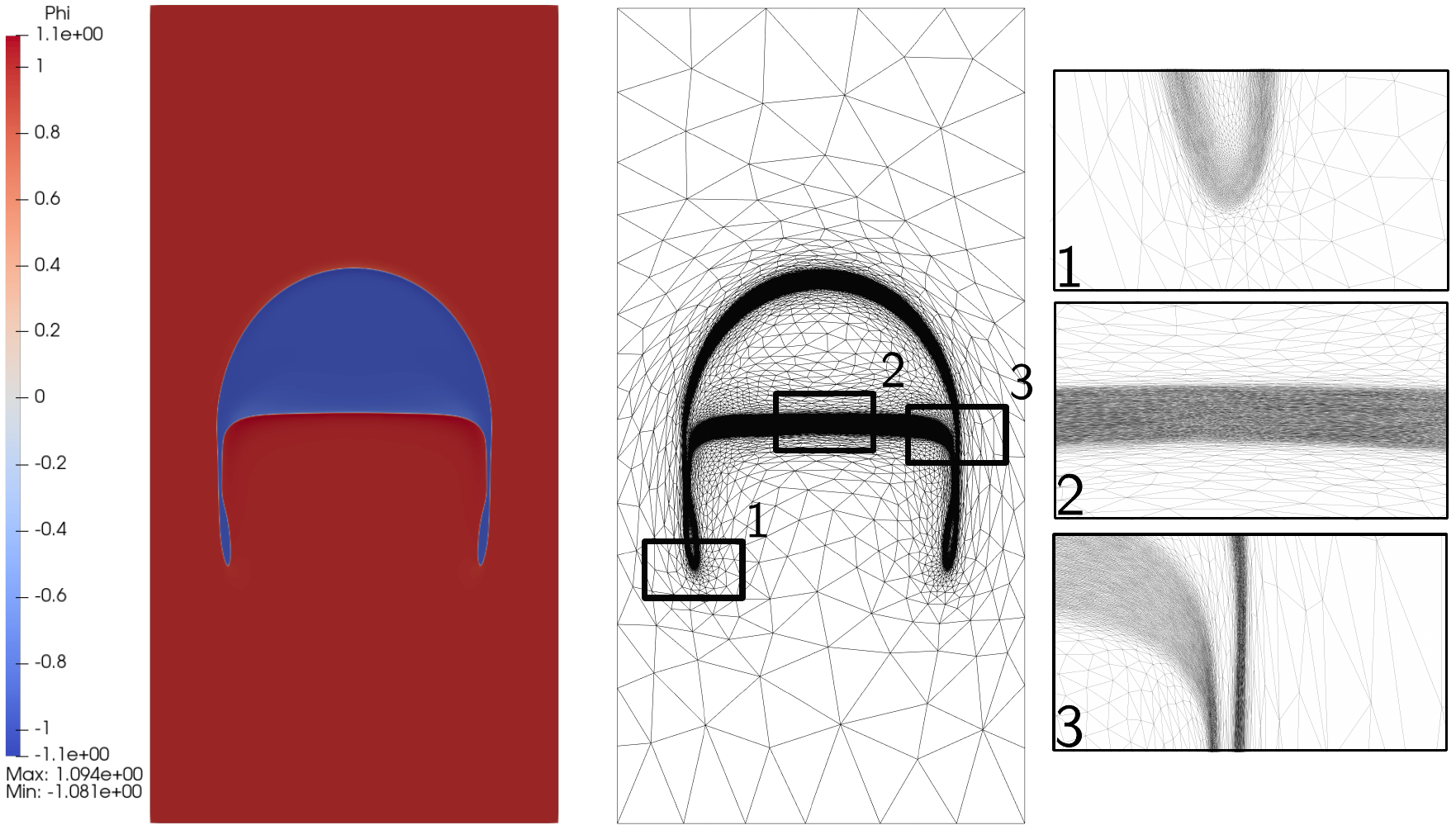}
  \caption{Bubble benchmark 2 - Phi and last adapted mesh for 2048k vertices (give associated space-time complexity). Final mesh has X vertices and x triangles.}
  \label{fig:bubble2_phiMesh}
\end{figure}

\begin{figure}[h!]
  \centering
  \includegraphics[width=\linewidth]{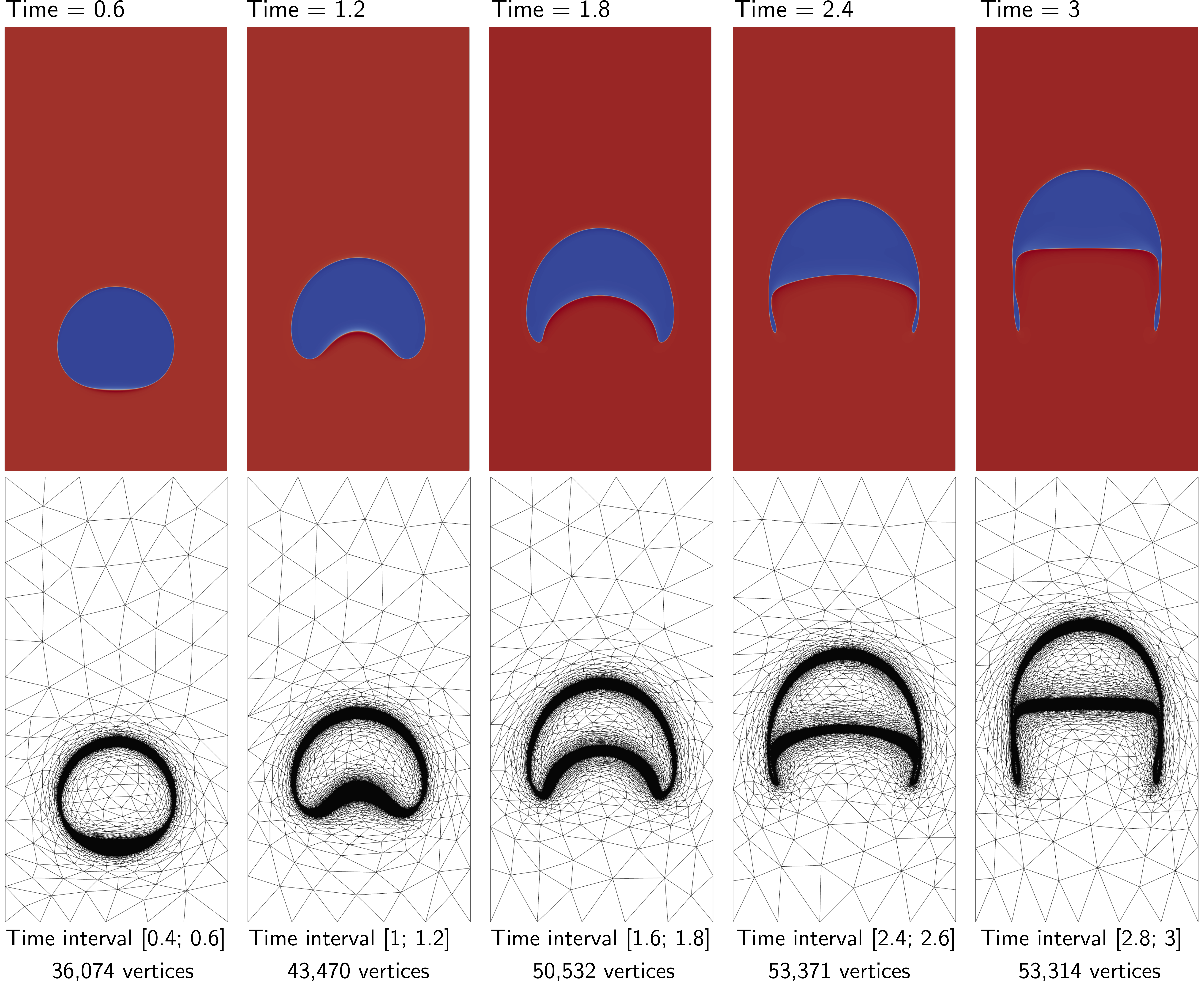}
  \caption{Bubble benchmark 2 - Phi and last adapted mesh for $\epsilon = 6.25 \times 10^{-4}$ and for 2048k vertices (give associated space-time complexity).}
  \label{fig:bubble2_rising}
\end{figure}

\paragraph{Benchmark 2}
The phase marker and final mesh for the second benchmark are shown on \Cref{fig:bubble2_phiMesh},
with an extra refinement step reaching $\epsilon = 0.000625$.
Here the over- and undershoot are more pronounced, with $\phi(T) \in [-1.081, 1.094]$.
The filaments stay attached to the rest of the bubble despite the thin interface.
The bubble rising at different times is depicted on \Cref{fig:bubble2_rising},
along with the adapted mesh for those time intervals. Each depicted solution is associated to the final instant
of each interval.

\begin{figure}[h!]
  \centering
  \includegraphics[width=0.8\linewidth]{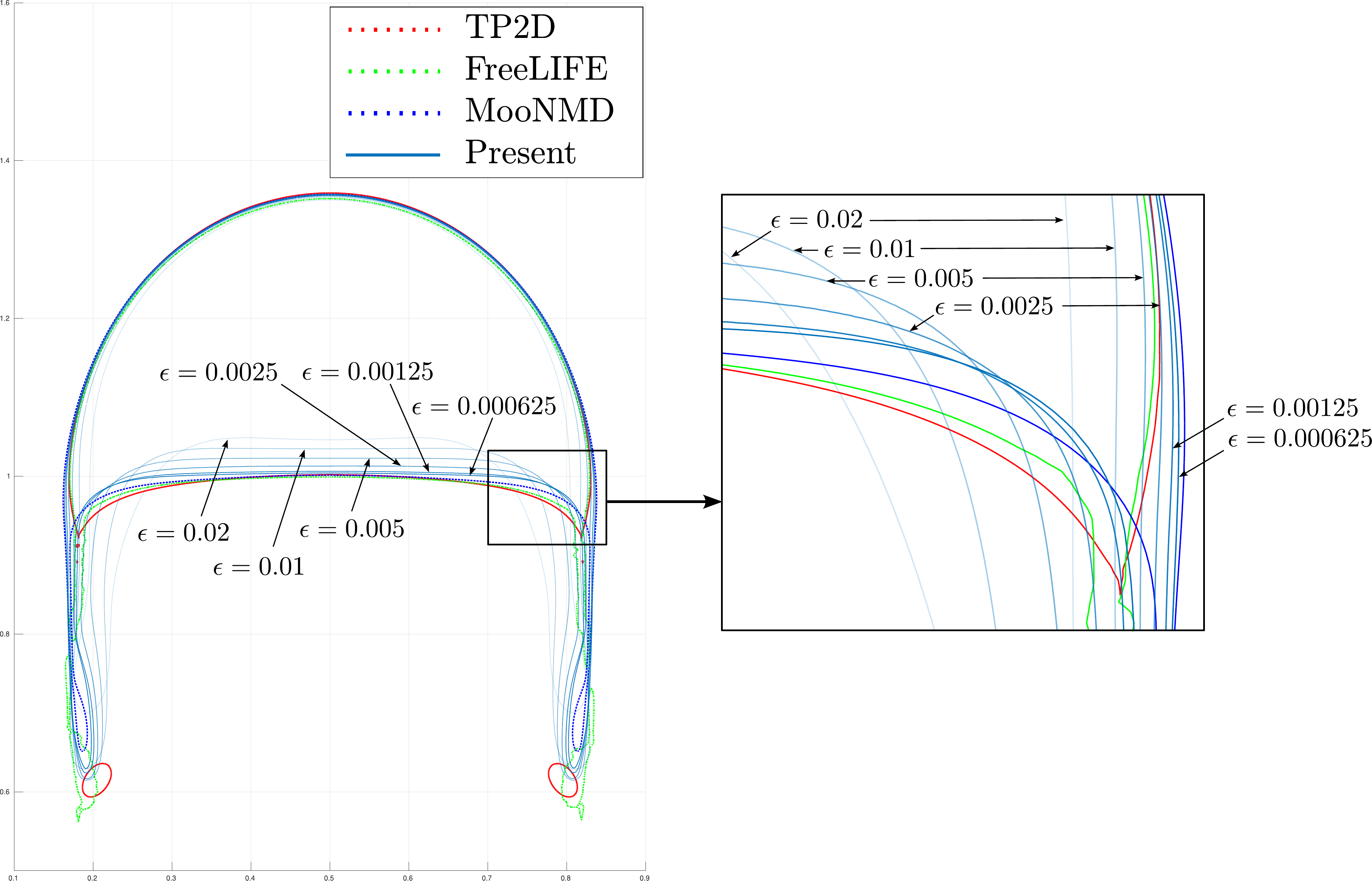}
  \caption{Bubble benchmark 2 - Comparison with the literature: shape of the bubble.}
  \label{fig:bubble2_shape}
\end{figure}

Comparison with the literature for the bubble shape and quantitative metrics are provided on \Cref{fig:bubble2_shape}
and right of \Cref{fig:bubblesIndicators} for decreasing $\epsilon$.
Once again, excellent agreement is obtained at finer interface thickness.

We should point out that on a uniform grid, a naive run for the smallest interface thickness would require
a mesh size $h = 1/1024$ if taking $\epsilon = 0.64 h$ as in \citep{brunk2025simple},
yielding a uniform mesh with $1024 \times 2048 = 2,097,152$ vertices.
Here, the finest anisotropic mesh for this $\epsilon$ contains 53,314 vertices,
thus a reduction by a factor of almost 40, and even more in terms of degrees of freedom.



\subsubsection{Convergence study for benchmark 1}
Lastly, we assess the convergence of the mesh adaptive process for the first benchmark
with interface thickness fixed to $\epsilon = \epsilon_0 = 0.02$.
As in \Cref{sec:mms_adaptive}, two convergence studies are performed,
one at fixed number of time intervals $\nI = 15$ and increasing mesh density
$\ncalavg = 64$k $\to 4.1$M, and one at fixed mesh density $\ncalavg = 128$k
and increasing number of intervals $\nI = 2 \to 128$.
The time step is $\Delta t = T/750$ for the first study and $T/1280$ for the second.
As no exact solution is available, the reference solution is the numerical solution
on the finest sequence of meshes for each study.
Unlike for the manufactured solution, here convergence is nearly identical
for both configurations, likely due to the accumulation of transfer error
on the four variables which penalizes the use of more intervals,
and we observe a rate of convergence of at least 3 for
the phase marker $(\sim \nst^{-{3/3}})$.
We also note that even though the meshes are adapted for $\phi$ only, a
convergence at order 2+ is obtained for the other scalar fields.

\begin{figure}[h!]
  \centering
  \includegraphics[width=\linewidth]{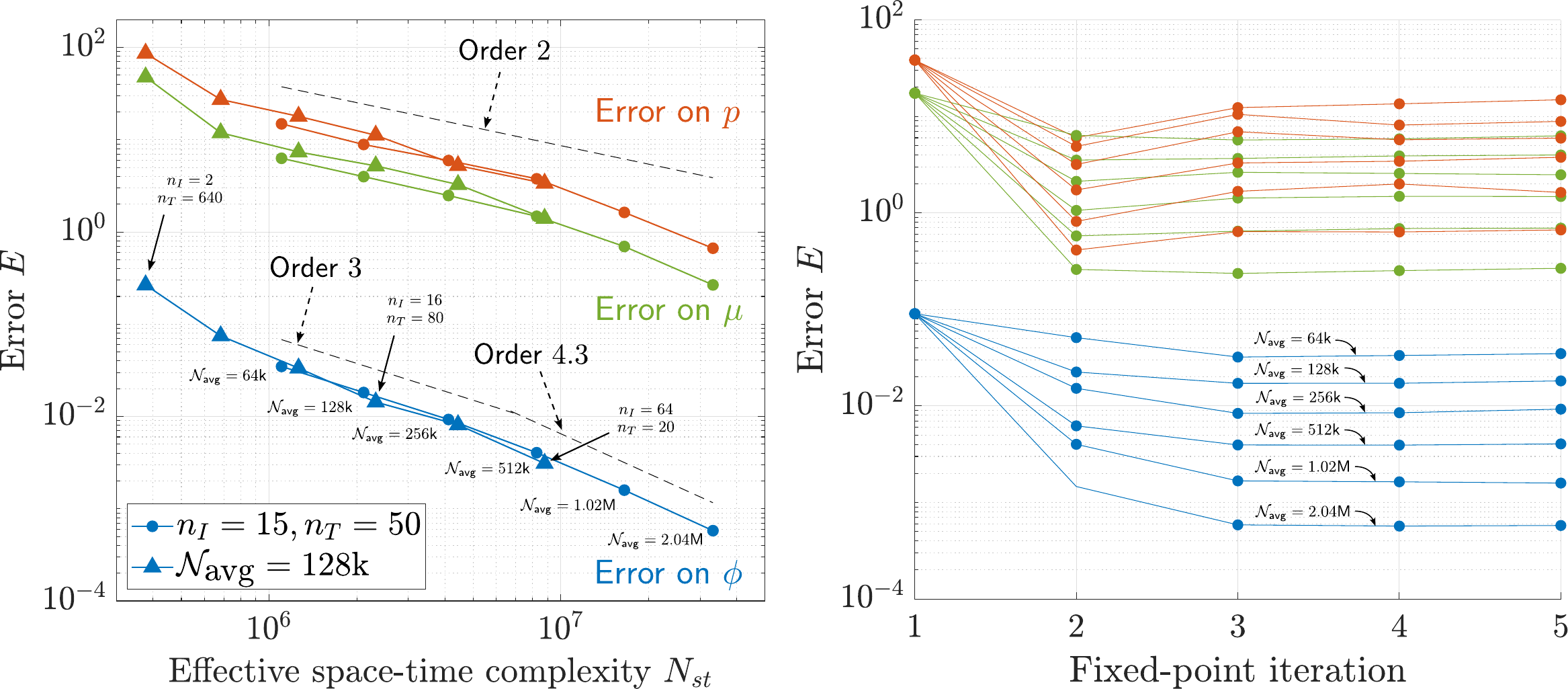}
  \caption{Bubble benchmark 1. Left: Convergence study at constant time step on $p, \phi$ and $\mu$: fixed number of intervals and increasing mesh density (circles) vs. fixed mesh density and increasing number of intervals (triangles). Right: Evolution of the errors with respect to the finest mesh solution vs. number of fixed-point iterations.}
  \label{fig:bubble1_convergence}
\end{figure}

The evolution of the errors against the number of fixed-point iterations
is given on the right of \Cref{fig:bubble1_convergence} for the study at constant $\nI$.
It points out a fast convergence of the global fixed-point method, as the errors
stabilize after 2 to 3 fixed-point iterations (that is, after the complete simulation has been run 2 to 3 times).

\section{Conclusion}

We presented a fully automatic mesh adaptive method for unsteady multiphase simulations
with the Cahn-Hilliard Navier-Stokes framework, based on the existing global transient
fixed-point method.
By running the complete multiphase simulation a few times,
the mesh is adapted to the unsteady features of the phase marker and
captures the fluid-fluid interface without lag, nor needing a priori mesh refinement.

In the context of a diffuse interface model,
anisotropic mesh adaptation allows for both an accurate representation of the
fluid-fluid interface, and a dynamic refinement of the interface
thickness. The initial simulation can be computed on a coarse grid with a large
$\epsilon$, then converge to a finely resolved interface in a completely
automatic way, making it particularly appealing.

The procedure was verified with manufactured solutions and applied to the rising bubbles benchmark.
Results at low epsilon compare well to front-tracking results
and only use a fraction of the mesh elements requires with a uniform grid.

Regarding the discretization and the mesh adaptive procedure, various improvements can be mentioned as directions of future research:
\begin{itemize}
\item using a conservative solution transfer, to limit the accumulation of transfer error;
  \item going to higher-order numerical scheme enjoying faster convergence, e.g., quadratic or higher discretization of the phase marker;
  \item applying the time and space adaptive procedure described in \citep{sauvage2024space}
  for a fully adaptive numerical scheme;
  \item considering high-order, metric-conforming mesh adaptive techniques (i.e., curved meshes), for a better representation of the resolved fields.
\end{itemize}

%


 \bibliographystyle{elsarticle-num-names}
 \bibliography{references}



%
%
%
\end{document}